\documentclass[prl,onecolumn,showpacs,twoside,showpacs]{revtex4}

\usepackage{graphicx}
\usepackage{dcolumn}
\usepackage{bm}
\usepackage{color}
\usepackage{amsmath}
\usepackage{cases}
\usepackage{ulem}

\begin{document}\large
%\begin{CJK*}{GBK}{song}

\title{Efficient Matching Boundary Conditions of Two-dimensional Honeycomb Lattice for Atomic Simulations}
\author{Baiyili Liu$^1$, Songsong Ji$^2$, Gang Pang$^3$, Shaoqiang Tang$^2$ and Lei Zhang$^{4,5}$}
\email{zhanglei@imech.ac.cn}
\affiliation{$^1$ College of Physics and Electronic Engineering, Center for Computational Sciences, Sichuan Normal University, Chengdu 610066, China}
\affiliation{$^2$ HEDPS, CAPT, and LTCS, College of Engineering, Peking University, Beijing 100871, China.}
\affiliation{$^3$ School of Mathematical Science, Beihang University, Beijing 102206, China}
\affiliation{$^4$ The State Key Laboratory of Nonlinear Mechanics, Institute of Mechanics, Chinese Academy of Sciences, Beijing 100190, China}
\affiliation{$^5$ School of Engineering Science, University of Chinese Academy of Sciences, Beijing 100049, China}

\date{\today}
\begin{abstract}
In this paper, we design a series of matching boundary conditions for a two-dimensional compound honeycomb lattice, which has an explicit and simple form, high computing efficiency and good effectiveness of suppressing boundary reflections. First, we formulate the dynamic equations and calculate the dispersion relation for the harmonic honeycomb lattice, then symmetrically choose specific atoms near the boundary to design different forms of matching boundary conditions. The boundary coefficients are determined by matching a residual function at some selected wavenumbers. Several atomic simulations are performed to test the effectiveness of matching boundary conditions in the example of a harmonic honeycomb lattice and a nonlinear honeycomb lattice with the FPU-$\beta$ potential. Numerical results illustrate that low-order matching boundary conditions mainly treat long waves, while the high-order matching boundary conditions can efficiently suppress short waves and long waves simultaneously. Decaying kinetic energy curves indicate the stability of matching boundary conditions in numerical simulations.
\end{abstract}

\maketitle
\date{today}
%%%%%%%%%%%%%%%%%%%%%%%%%%%%%%%%%%%%%%%

\section{1. Introduction}
Molecular dynamics method becomes an instrumental tool for researching physical properties of materials at nano-scale \cite{S. Hajilar2015,A. Kumar2020}, or simulating mechanical behaviours of a crystal lattice like moving dislocation  \cite{S. Groh2009,E. Oren2016}, vacancy migration \cite{Y. Lee2011,G. Berdiyorov,M. Noshin}, or crack propagation \cite{C. Rountree2002,P. Budarapu2015}. Although the rapid development of computing power, full atomic simulations for mechanical materials or devices in macroscopic level still costs a great deal of computing time and large computing internal storages. To improve computing efficiency, people usually take multiscale computational methods by separating the whole region into a molecular dynamics subdomain and a macroscopic subdomain. The core domains that containing interesting physical problems are simulated by atomic simulations, and other domains are described by continuum model \cite{Kadowaki2005, R. Miller2009, V. Yamakov2014}, or coarse grid equations \cite{S. Tang2006}. For finite atomic simulation region, appropriate artificial boundary conditions need to be proposed to make dynamic equations closed. Different from the real physical boundary, the artificial boundary condition is more complicated and needs to be specially designed for different crystal lattices.  Boundary reflections will occur and disturb the accuracy of simulations, if they are not be treated properly. What's more, inappropriate artificial boundary conditions may bring numerical instability making the simulation impossible.

In recent years, many artificial boundary conditions have been proposed for discrete lattices, such as time history kernel treatment (THK) \cite{Adelman&Doll 1974, Adelman&Doll 1976}, perfectly matched layer (PML) \cite{Berenger1994}, variational boundary condition (VBC) \cite{W. E2001, W. E2002}, velocity interficial condition (VIC) \cite{S. Tang2008}, almost exact method (ALEX) \cite{G. Pang ALEX 2012} and matching boundary condition (MBC) \cite{X. Wang2010, X. Wang2013}.

The THK is a kind of exact artificial boundary condition written in a convolutional form over time \cite{Adelman&Doll 1974, Adelman&Doll 1976}. The convolution kernel function is obtained analytically for specific lattices \cite{W. Liu2000, W. Liu2003, W. Liu2004, M. Dreher2008, G. Pang ALEX 2011}. Although it is an exact boundary condition, the calculation of the convolution requires increasing memory and computations with time, which means an immense cost for long time simulations.

To improve efficiency, people developed approximated ones, which are local in both space and time. The PML method introduces an absorbing layer with dissipative terms to absorb electromagnetic waves with various frequencies in all directions, first proposed by Berenger in solving Maxwell equations \cite{Berenger1994}. It has been used in many problems, such as the elastic waves \cite{W. Chew1996, Hastings1996} and the nonlinear Schr\"{o}dinger equations \cite{C. Zheng2007}. Li et al. introduce the PML method into molecular dynamics of discrete crystal lattices \cite{S. Li2005, S. Li2006} . In the PML method, the choice of coefficients usually rely on experiences, which is difficult to achieve an accurate suppression of boundary reflections.

The VBC is a kind of local boundary condition represented by a linear combination of current and former moments of atomic displacements proposed by E and Huang \cite{W. E2001, W. E2002}.  Li and E extend the VBC to high-dimensional lattices \cite{X. Li2006, X. Li2007}. The undetermined coefficients of the VBC are calculated by optimizing reflection coefficients globally. Thus the VBC can well absorb short waves with high computing efficiency, but the procedure of optimization is thus trial and needs a large amount of calculations.

The VIC proposed by Tang is based on one-wave traveling wave equation \cite{S. Tang2008}. It is written in a linear combination of displacements of atoms near the boundaries. Except for linear lattices, the VIC is also applicable for simple lattices with weak nonlinearity.  Furthermore, Pang and Tang propose the ALEX method \cite{G. Pang ALEX 2011, G. Pang ALEX 2017}, which is one of the most exact artificial boundary conditions. The core technique of the ALEX is to calculate coefficients of the linear combination to approximate the convolution kernel function, which is based on the approximation theory of high order Bessel functions by low order ones \cite{G. Pang Bessel 2017}. The ALEX is highly accurate, but the derivation of approximating kernel function is very complicated, which limits the extension of the ALEX to high dimensional lattices.

The MBC is a kind of accurate artificial boundary condition proposed by Wang and Tang for low-dimensional simple lattices \cite{X. Wang2010, X. Wang2013}. The MBC adopts an explicit and simple form represented by a linear combination of displacements and velocities of atoms near the boundary. The coefficients are determined by matching dispersion relation at specific wave numbers without empirical parameters.
Compared with the VBC, the off-line calculation of coefficients is easy and fast. The MBC is also local in space and time, and thus has high computing efficiency remaining a high accuracy. In contrast to the ALEX, it is easy to extend the MBC to multiple-dimensional cases, such as the triangle lattice \cite{B. Liu2017}, the three-dimensional face centered cubic lattice \cite{M. Fang2012} and body centered cubic lattice \cite{M. Fang2013}. Ji and Tang verify the stability of the MBC for atomic chains \cite{S. Ji2014}. Tang and Liu further apply the MBC to design heat bath method for atomic chains \cite{B. Liu2015, B. Liu2020}. Zhang et.al. use the MBC to perform multiscale simulations at finite temperature \cite{L. Zhang2022}. In the present work, we extend the MBC to the honeycomb lattice, one kind of two-dimensional compound lattice.

In practical applications, many materials are two-dimensional or three-dimensional compound lattices rather than simple lattices. Thus designing an efficient artificial boundary condition for a compound lattice is very important. Among various compound lattices, the honeycomb lattice is typical and applicable for atomic simulations of real materials, such as graphene and carbon nanotube. Some artificial boundary conditions have been applied to the honeycomb lattice. For example, the THK is applied to multiple scale simulations of the honeycomb lattice \cite{G. Wagner2004}, or studying the indentation of graphene \cite{S. Medyanik2006}. The PML is improved for thin dispersive layers of graphene \cite{X. Yu2012}. Considering the advantages of the matching boundary condition, namely an explicit and simple form, high computing efficiency and good effectiveness of suppressing boundary reflections, we are going to systematically design the matching boundary conditions for a honeycomb lattice. It can be a basis for studying materials with complex compound lattice structures, and also be of great importance for the multiscale computation.

The rest of the paper is organized as follows. In section 2, we formulate dynamic equations and calculate the dispersion relation for the harmonic honeycomb lattice. Then, we design a series of matching boundary conditions and analyze reflection coefficients. Numerical simulations are performed, investigating the effectiveness of different orders of matching boundary  conditions  for the harmonic honeycomb lattice and the nonlinear case with FPU-$\beta$ potential in section 3. Some concluding remarks are presented in the last section. Detailed derivations of boundary coefficients are shown in Appendix.

\section{2. Lattice dynamics model and boundary condition}

In this section, we first describe the governing equations for the honeycomb lattice with nearest neighbouring harmonic interaction. Then we calculate the dispersion relation and analyze the group velocity in the extended Brillouin zone. In view of matching dispersion relation, we design four representative matching boundary conditions and analyze the reflection suppression effects by reflection coefficients.

\subsection{2.1  Dynamic equations and dispersion relation}
%\begin{figure}
%  \centering
%  \includegraphics[width=6cm]{triangle_atom1.eps}
%  \caption{Triangular lattice with nearest neighboring interaction.}\label{f:atom1}
%\end{figure}

\begin{figure}
  \centering
  \includegraphics[width=14cm]{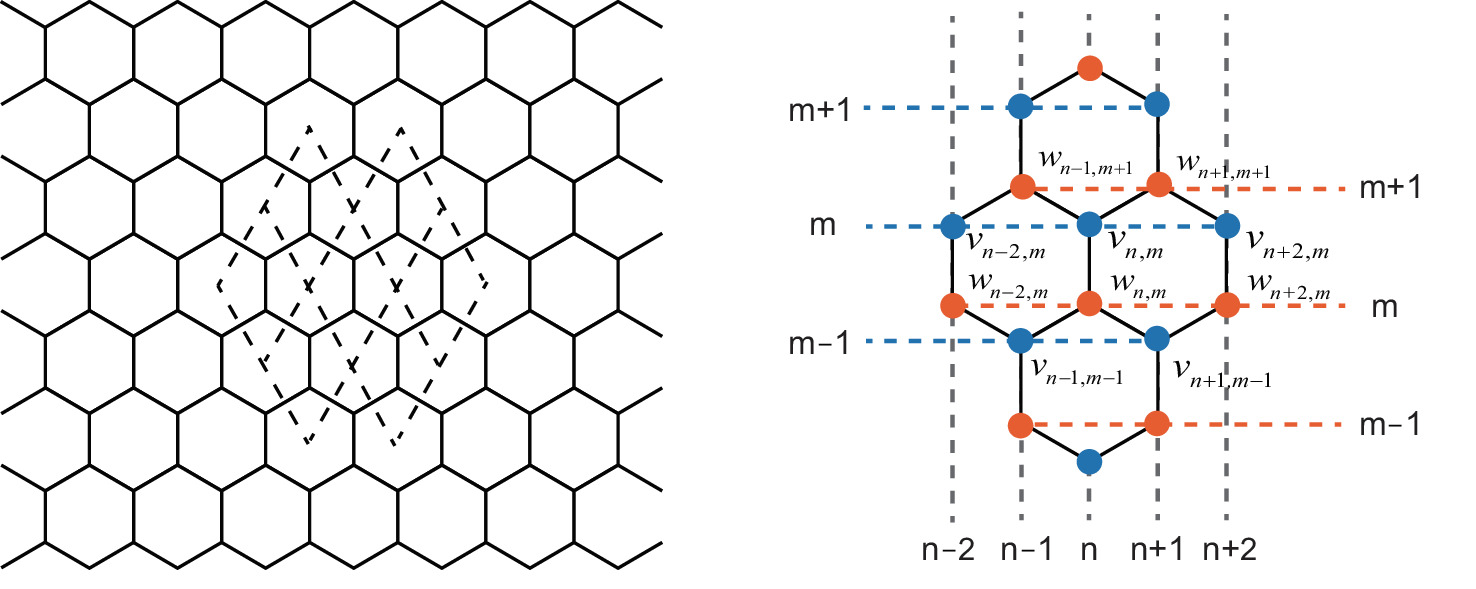}%
  \caption{Left: schematic plot of the honeycomb lattice. A diamond drawn by dotted lines represents a unit cell; Right: illustration for numbering atoms. The first class of atoms are represented by blue points with the displacements denoted by $v_{n,m}$. The second class of atoms are represented by red points with the displacements denoted by $w_{n,m}$. }
  \label{f:schematic}
\end{figure}

See the schematic plot of the honeycomb lattice in the left subplot of Fig. \ref{f:schematic}. A diamond drawn by dotted lines represents a unit cell. Each unit cell contains two classes of atoms, owning different symmetries and different dynamics, although they belong to the same kind of atoms, such as carbon atoms in the graphene. To better distinguish two classes of atoms, we use blue and red points for drawing in the right subplot. We also introduce two symbols $v_{n,m}$ and $w_{n,m}$ to represent their displacements respectively. The subscript $n$ in horizontal direction and $m$ in vertical direction are used to number atoms, as shown in the right subplot.

%Notice that there exist two classes of atoms in the lattice represented by blue and red points in the right subplot due to different positions and the corresponding different dynamics, although they belong to the same kind of atoms, such as carbon atoms in the graphene. That is to say, the honeycomb lattice is a two-dimensional compound lattice in spite of only one kind of atoms in the lattice. A diamond drawn by dotted lines represents a unit cell containing these two classes of atoms as shown in the left subplot. In this paper, we introduce two symbols $v_{n,m}$ and $w_{n,m}$ to represent their displacements respectively. The subscript $n$ in horizontal and $m$ in vertical are used to number atoms, as shown in the right subplot.

Consider out-of-plane motion \cite{S. Ji2017}, namely, all atoms move perpendicularly to the crystal plane. Thus the displacements and velocities of atoms  are scalars in the vertical direction. Motion of the system is governed by the Newton equations
\begin{eqnarray}\label{eq:motion}
  m \ddot{v}_{n,m}=\frac{\partial U}{\partial v_{n,m}}, ~~~~
  m \ddot{w}_{n,m}=\frac{\partial U}{\partial w_{n,m}}. \label{eq:harmonic_newton}
\end{eqnarray}
Here, $m$ is the mass of atoms, $U$ is a potential depending on the displacement difference for all interacting atoms. The variables $v_{n,m}$ and $w_{n,m}$ represent displacements of two classes of atoms away from their equilibrium positions respectively.

Consider each atom interacts with three nearest adjacent atoms under a linearized potential
\begin{eqnarray}\nonumber
&&U=\sum_{n,m}\frac{k}{2} \Big[ (w_{n-1,m+1}-v_{n,m})^2+(w_{n+1,m+1}-v_{n,m})^2
+(w_{n,m}-v_{n,m})^2\\
&&~~~~~~+(w_{n,m}-v_{n-1,m-1})^2+(w_{n,m}-v_{n+1,m-1})^2\Big]. \label{eq:harmonic_potential}
\end{eqnarray}
Here, $k$ is the linear force constant. Substitute the linearized potential (\ref{eq:harmonic_potential}) into the Newton equations (\ref{eq:harmonic_newton}), and then rescale time by $\sqrt{m/k}$ and the displacement by lattice constant $a$. The corresponding non-dimensional Newton motion equations read
\begin{eqnarray}\label{eq:harmonic_newton1}
   \ddot{v}_{n,m}&=&w_{n-1,m+1}+w_{n+1,m+1}+w_{n,m}-3v_{n,m}, \\\label{eq:harmonic_newton2}
   \ddot{w}_{n,m}&=&v_{n-1,m-1}+v_{n+1,m-1}+v_{n,m}-3w_{n,m}.
\end{eqnarray}

%Take the commonly used Morse potential for example, the linear force constant of monocrystalline silicon is equal to $k=2\varepsilon \alpha^2=*$ with parameter $\varepsilon=*$ and $\alpha=*$.
\begin{figure}
 \centering
  \includegraphics[width=9cm]{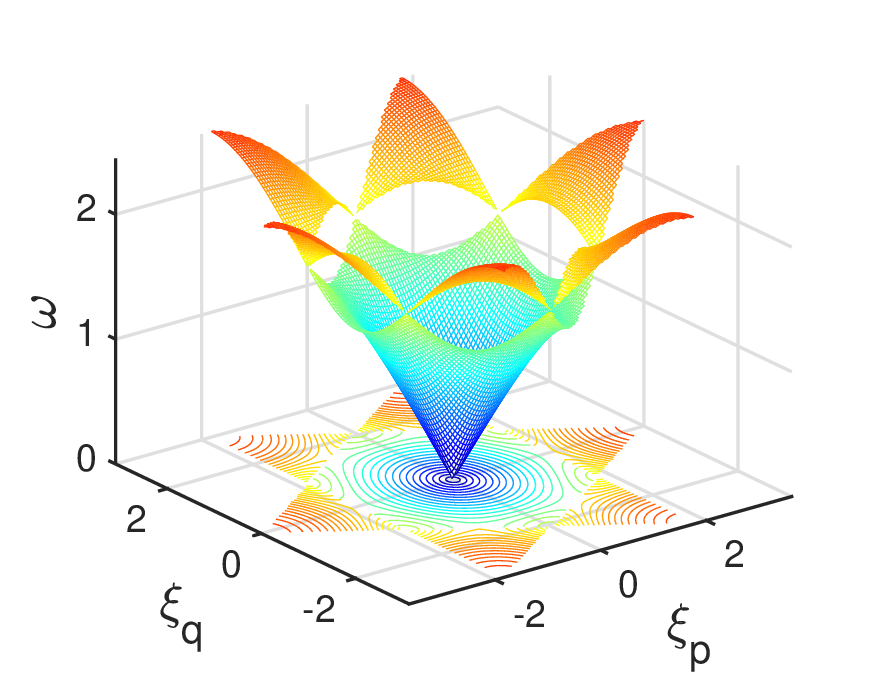}
  \includegraphics[width=8cm]{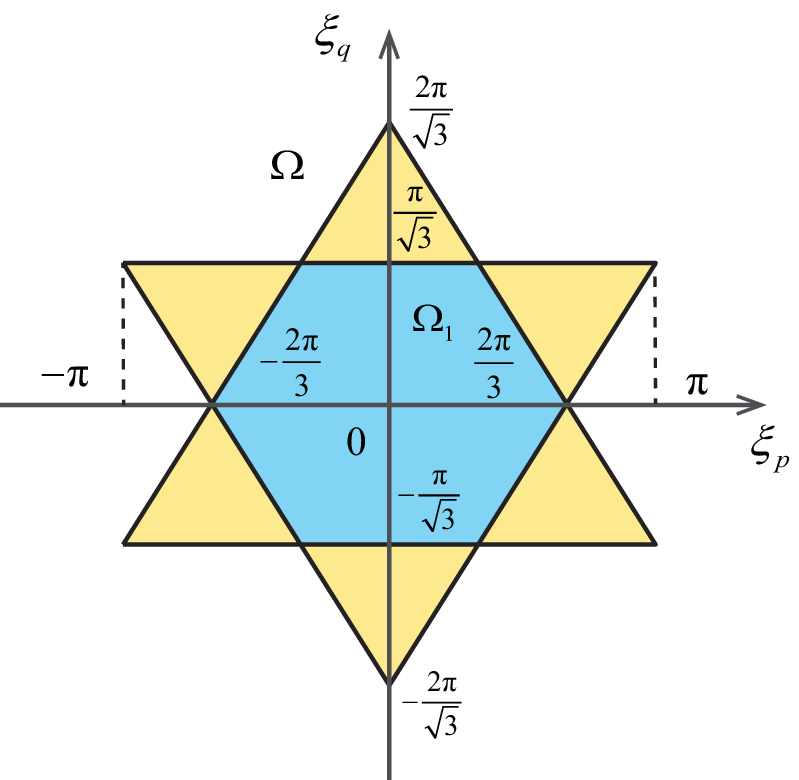}
  \caption{Left: dispersion relation of the honeycomb lattice; Right: the blue domain is the reduced Brillouin zone $\Omega_1$, the whole domain is the extended Brillouin zone $\Omega$, and the yelllow domain corresponds to the optical branch $\Omega\setminus\Omega_1$.}\label{f:dispersion}
\end{figure}

Substituting a monochromatic wave with a wave vector $(\xi_p,\xi_q)$ in the form of $v_{n,m}=A_1 e^{i(\omega t + \xi_pn + \sqrt{3}\xi_qm)}$ and $w_{n,m}=A_2 e^{i(\omega t + \xi_pn + \sqrt{3}\xi_qm)}$ into the Newton equations (\ref{eq:harmonic_newton1}) and (\ref{eq:harmonic_newton2}), we obtain the eigen equations
\begin{eqnarray}
\left(
\begin{array}{cc}
\omega^2-3 & 1+2\cos \xi_p e^{i\sqrt{3}\xi_q} \\
1+2\cos \xi_p e^{-i\sqrt{3}\xi_q} &  \omega^2-3
\end{array}
\right)
\left(
\begin{array}{c}
A_1 \\
A_2
\end{array}
\right)=0. \label{eq:A1A2}
\end{eqnarray}
Here, $A_1$ and $A_2$ are the amplitudes of lattice waves.

To guarantee the existence of  nontrivial solution of Eq.(\ref{eq:A1A2}), it requires that the determinant equals zero
\begin{eqnarray}
\left|
\begin{array}{cc}
\omega^2-3 & 1+2\cos \xi_p e^{i\sqrt{3}\xi_q} \\
1+2\cos \xi_p e^{-i\sqrt{3}\xi_q} &  \omega^2-3
\end{array}
\right|
=0. \label{eq:det}
\end{eqnarray}
This leads to the dispersion relation
\begin{eqnarray}
\omega_{p,q,s}^2=3+(-1)^s\sqrt{1+4\cos \xi_p \cos\sqrt{3}\xi_q + 4\cos^2 \xi_p}. \label{eq:dispersion}
\end{eqnarray}

The dispersion relation has two branches: $s=1$ and $s=2$ correspond to the acoustic branch and the optical branch respectively. The curved surface of the dispersion relation is shown in the left subplot of Fig. \ref{f:dispersion}. The lower part of image is the acoustic branch. When the wave vector $\pmb{\xi} \rightarrow 0$, the frequency of acoustic branch equals zero. Due to periodicity in the dispersion relation, a common choice is the extended Brillouin zone $\Omega$, which is a shape of hexagram as shown in the right subplot of Fig. \ref{f:dispersion}. In the middle region of $\Omega$, it is the reduced Brillouin zone $\Omega_1$  defined by
\begin{eqnarray}
\Omega_1=\Big\{ (\xi_p,\xi_q) \Big| \sqrt{3}|\xi_p|+|\xi_q|\leq \frac{2\pi}{\sqrt{3}} \cap |\xi_q|\leq \frac{\pi}{\sqrt{3}}  \Big\}.\label{eq:Brillouin}
\end{eqnarray}
The rest domain $\Omega\setminus\Omega_1$ corresponds to the optical branch. Substituting the dispersion relation (\ref{eq:dispersion}) into the eigen equations (\ref{eq:A1A2}), we can obtain the proportion of amplitudes
\begin{eqnarray}
\left(
\begin{array}{c}
A_1 \\
A_2
\end{array}
\right)
=
\left(
\begin{array}{c}
1+2\cos \xi_p e^{i\sqrt{3}\xi_q} \\
- \omega_{p,q,s}^2 +3
\end{array}
\right). \label{eq:amplitude}
\end{eqnarray}

\begin{figure}
 \centering
  \includegraphics[width=8cm]{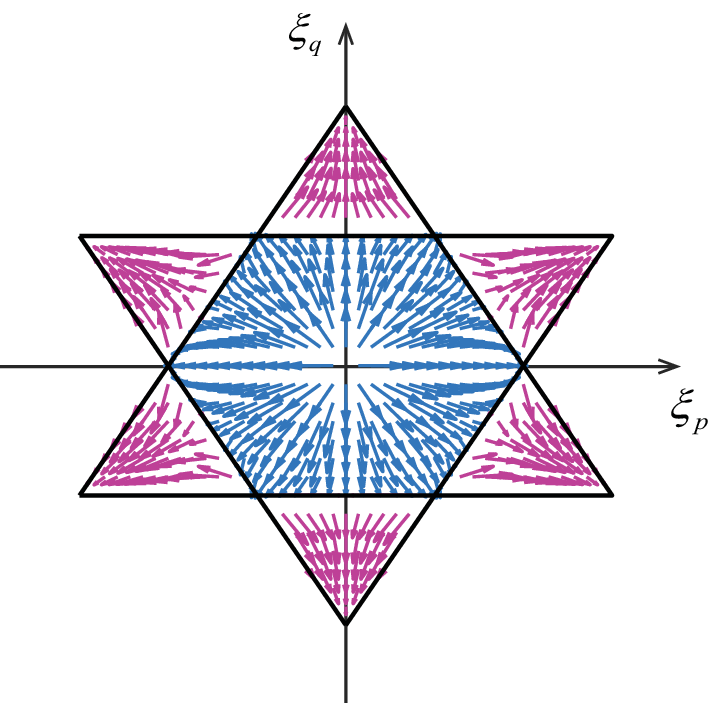}
  \includegraphics[width=8.5cm]{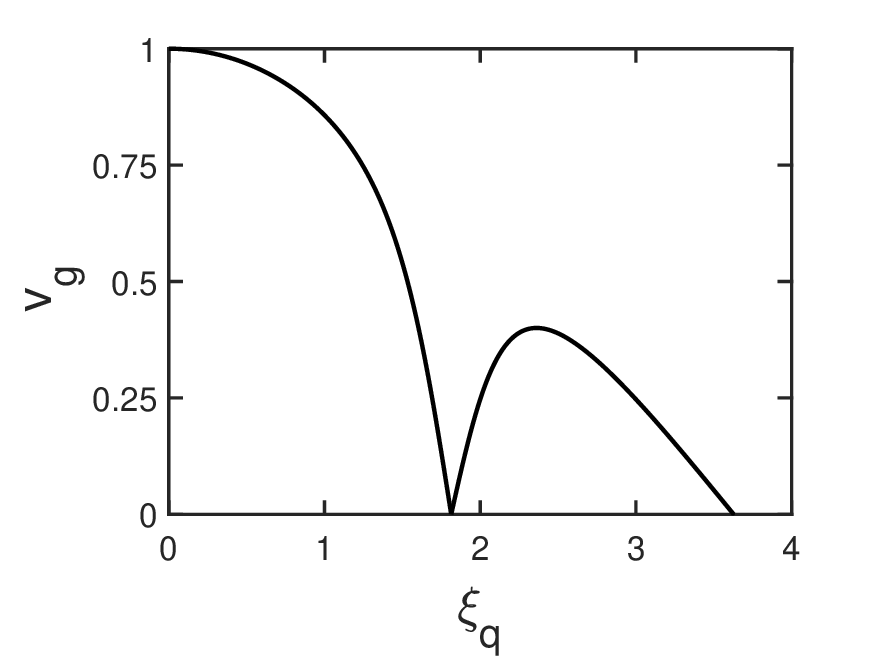}
  \caption{Left: vector graphic of the group velocity; Right: the group velocity along the vertical direction with $\xi_p=0$.}\label{f:group_velocity}
\end{figure}

Then we calculate the group velocity vector for giving a guidance to design the matching boundary conditions. Taking the derivatives of the dispersion relation function, the group velocity vector is written as
\begin{equation}
\nabla_{(\xi_p,\xi_q)}\omega_{p,q}=\left(\frac{8 \cos \xi_p \sin \xi_p + 4 \cos \sqrt{3} \xi_q \sin \xi_p}{4\omega_{p,q,s}  (3-\omega_{p,q,s}^2)},\frac{\sqrt{3}\sin \sqrt{3} \xi_q \cos \xi_p}{4\omega_{p,q,s} (3-\omega_{p,q,s}^2)}\right).
\end{equation}

The vector graphic of the group velocity is displayed in the left subplot of Figure. \ref{f:group_velocity}. The arrow represents the direction of wave propagation. It can be seen that the length of the arrow is longest around $\pmb{\xi}$ $\rightarrow 0$, which illustrates the group velocity is the biggest at the long wave limit. Let the wave number $\xi_p=0$, we plot the group velocity along the vertical direction in the right subplot of Figure. \ref{f:group_velocity}. It can be observed that the group velocity has two local maximum values at $\xi_q=0$ and $\xi_q=2.4$. These critical points may be used for designing matching boundary conditions in the next subsection.

%The absolute value $|A_1|$ is equal to $A_2$. Thus in a unit cell, two kinds of atoms vibrate with the same amplitude, only phase difference between them.
%
%Substituting the amplitudes (\ref{eq:amplitude}) into Eq.(\ref{eq:vw}), the real solutions of lattice wave are written as
%\begin{eqnarray}
%&&v_{n,m}=|A_1| \cos(\omega_{p,q,s} t + \xi_pn + \sqrt{3}\xi_qm + \varphi), \\
%&&w_{n,m}=A_2 \cos(\omega_{p,q,s} t + \xi_pn + \sqrt{3}\xi_qm),
%\end{eqnarray}
%with phase $\varphi=\arctan( \frac{\text{Im} A_1 }{ \text{Re} A_1})=\frac{2\cos p \sin \sqrt{3} q}{1+2 \cos p \cos \sqrt{3} q}$.

\subsection{2.2 Design matching boundary conditions}
\begin{figure}
  \centering
  \includegraphics[width=14cm]{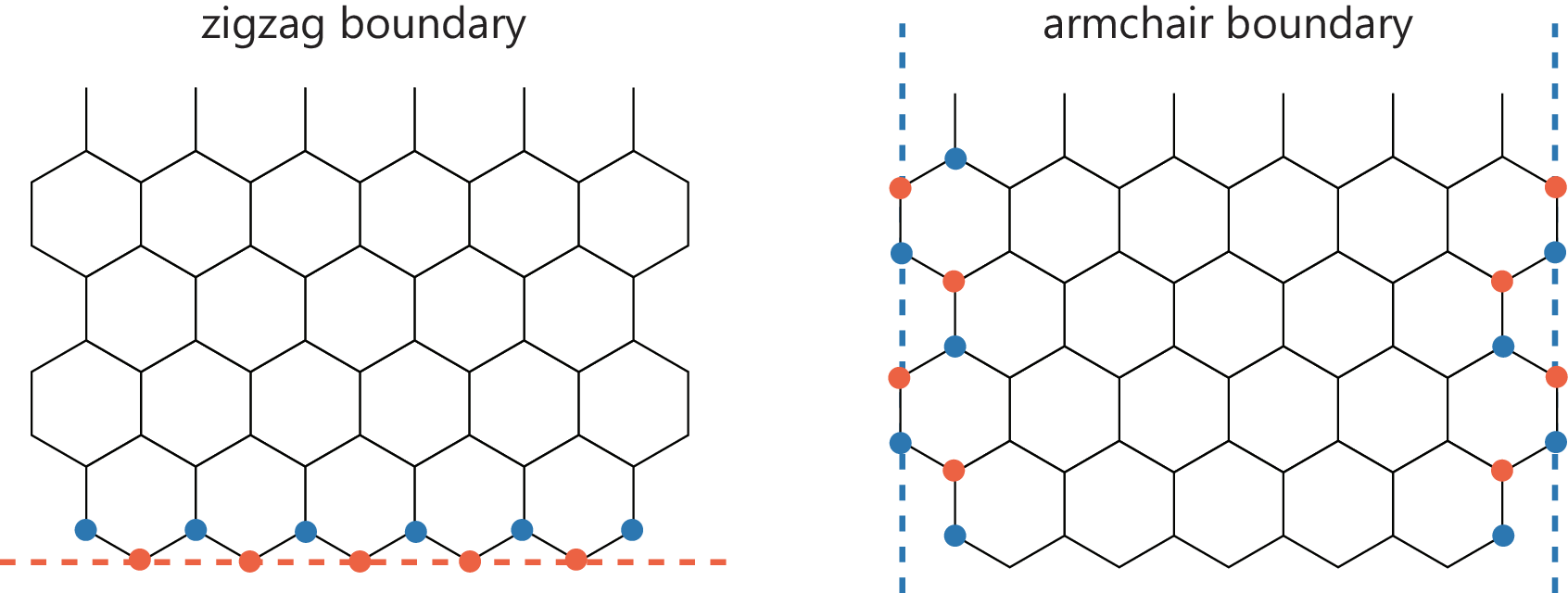}
  \caption{Schematic plot of two types of boundaries. Left: the zigzag boundary; Right: the armchair boundary.}\label{f:zigzag-armchair}
\end{figure}

\begin{figure}
  \centering
  \includegraphics[width=14cm]{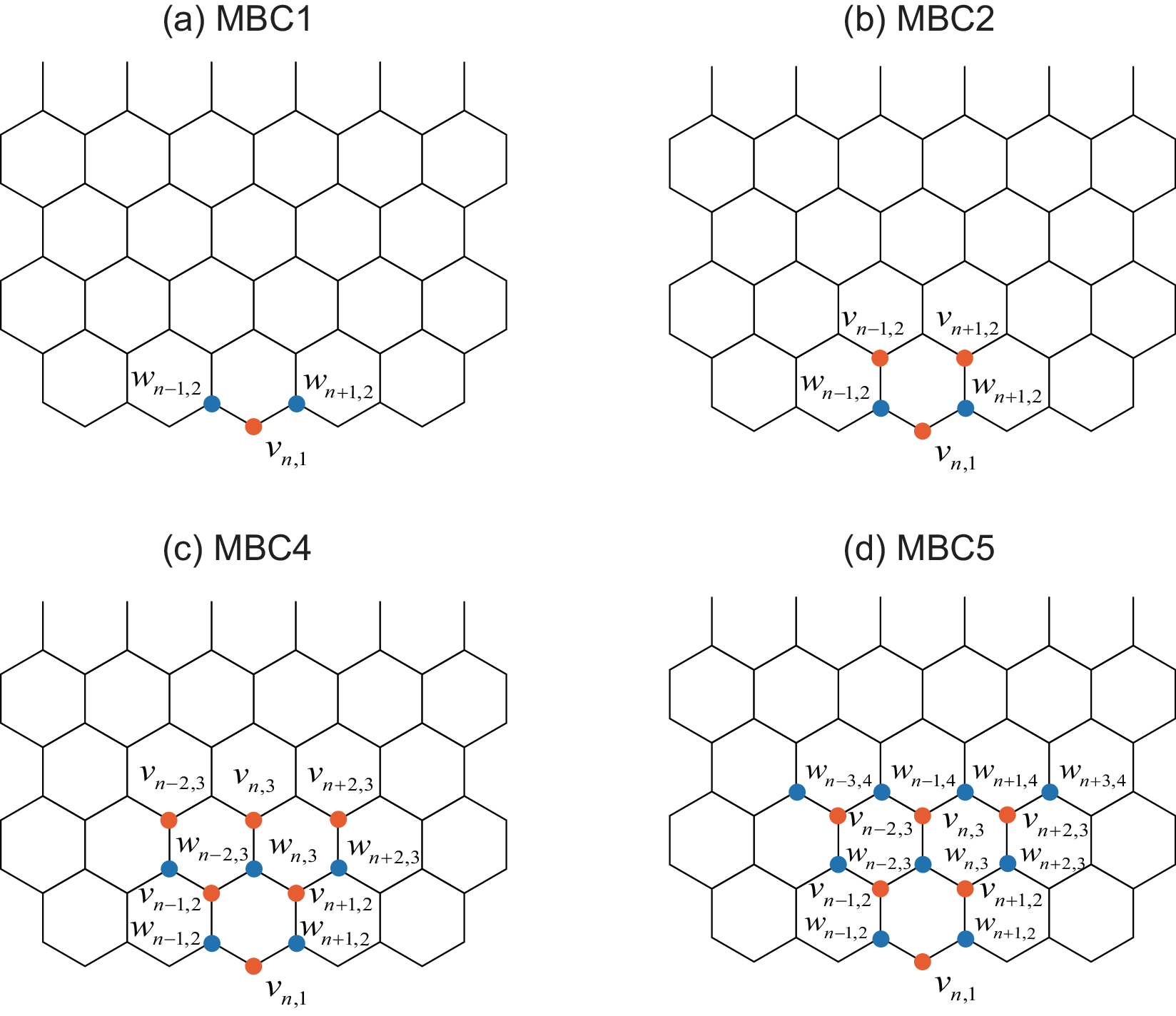}
  \caption{Schematic plot for the bottom matching boundary conditions: (a) MBC1; (b) MBC2; (c) MBC4; (d) MBC5.}\label{f:boundary-4}
\end{figure}

The honeycomb lattice has two types of boundaries \cite{S. Ji2017}. One is the zigzag boundary and another one is the armchair boundary as illustrated in Fig. \ref{f:zigzag-armchair}. The atoms near the zigzag boundary are arranged symmetrically along normal direction, which brings great convenience to design the matching boundary conditions. Thus, in this paper, we focus on designing four representative matching boundary conditions for the zigzag boundary as shown in Fig. \ref{f:boundary-4}: low-order MBC1 and MBC2 mainly for long waves in the acoustic branch, and high-order MBC4 and MBC5 for eliminating reflections of incident waves with broad band wave numbers and multi oblique angles. The number 'r' in the MBCr represents the number of layers involved in the matching boundary conditions.

Taking the bottom boundary as an example, the involved boundary atoms are shown by the solid circles in Fig. \ref{f:boundary-4}.  Let their displacements and velocities satisfy the following linear relation
\begin{eqnarray}\nonumber
b_{0,0}\dot{v}_{n,1}+\sum_{n',m'}b_{n'-n,2m'-3}\dot{w}_{n',m'}
+\sum_{n',m'}b_{n'-n,2m'-2}\dot{v}_{n',m'} \\
=\sum_{n',m'}c_{n'-n,2m'-3}w_{n',m'}+\sum_{n',m'}c_{n'-n,2m'-2}v_{n',m'} \Big. \label{eq:formula_mbc}
\end{eqnarray}
with subscripts
\begin{eqnarray}
n'=n-3, n-2, \cdots, n+2, n+3, ~~~~m'=1,2,3,4,
\end{eqnarray}
for the MBC5, which represent the numbering of resolving boundary atoms. We introduce $i=n'-n$ and $j=2m'-3$ as the subscript to number undetermined boundary coefficients $b$ and $c$. The integer $i$ is taken from $-3$ to $3$, and $j$ is taken from $0$ to $5$ in the MBC5. Considering the symmetry along normal direction, we take the coefficients $b_{i,j}=b_{-i,j}$ and $c_{i,j}=c_{-i,j}$. Let $b_{0,0}=1$, the other variables $b_{i,j}$ and $c_{i,j}$  are solved  by matching dispersion relation at  specific wave numbers.

The MBC$r$ ($r=1,2,3,4$) takes the first $r$ layer of atoms as shown in Fig. \ref{f:boundary-4}. For example, the MBC1 involves boundary atoms $v_{n,1}$, $w_{n-1,2}$ and $w_{n+1,2}$, accordingly $n'=n-1,n,n+1$ and $m'=1,2$, which results in three undetermined coefficients $b_{1,1}=b_{-1,1}, c_{0,0}$ and $c_{1,1}=c_{-1,1}$ ($b_{-1,1}=b_{1,1}$, $c_{-1,1}=c_{1,1}$ for symmetry). Thus the MBC1 takes form as
\begin{eqnarray}
\dot{v}_{n,1}+b_{1,1}\big(\dot{w}_{n-1,2}+\dot{w}_{n+1,2}\big)=
c_{0,0}v_{n,1}+c_{1,1}\big(w_{n-1,2}+w_{n+1,2}\big).
\label{eq:formula_mbc1}
\end{eqnarray}

Substituting the wave forms $v_{n,m} \sim A_1 e^{i(\omega t + \xi_p n + \sqrt{3}\xi_q m)}$ and $w_{n,m} \sim A_2 e^{i(\omega t + \xi_p n + \sqrt{3}\xi_q m)}$ into Eq.(\ref{eq:formula_mbc}), then discarding a common factor $e^{i(\omega t + \xi_pn + \sqrt{3}\xi_q )}$, we define a matching residual function

\begin{normalsize}
\begin{eqnarray}\nonumber
\displaystyle \Delta(\xi_p,\xi_q)&=&i \omega\bigg( A_1+A_2\sum_{n',m'}b_{n'-n,2m'-3}e^{i\big(\xi_p(n'-n)+\sqrt{3}\xi_q(m'-1)\big)} +A_1\sum_{n',m'}b_{n'-n,2m'-2}e^{i\big(\xi_p(n'-n)+\sqrt{3}\xi_q(m'-1)\big)}\bigg)\\
&&\displaystyle -\bigg(A_2\sum_{n',m'}c_{n'-n,2m'-3}e^{i\big(\xi_p(n'-n)+\sqrt{3}\xi_q(m'-1)\big)}+
A_1\sum_{n',m'}c_{n'-n,2m'-2}e^{i\big(\xi_p(n'-n)+\sqrt{3}\xi_q(m'-1)\big)}\bigg),
\end{eqnarray}
\end{normalsize}
which measures the performance of matching boundary conditions at specific wave vectors $(\xi_p, \xi_q)$.

As long waves usually occupy a major part of energy and have a relatively high speed, we preferentially eliminate boundary reflections at long wave limit. Let $\xi_p=0$, the matching residual function is written as
\begin{normalsize}
\begin{eqnarray}\nonumber
\displaystyle \Delta(0,\xi_q)&=&i \omega(0,\xi_q)\bigg( A_1(0,\xi_q)+A_2(0,\xi_q)\sum_{n',m'}b_{n'-n,2m'-3}e^{i\sqrt{3}\xi_q(m'-1)}
\displaystyle +A_1(0,\xi_q)\sum_{n',m'}b_{n'-n,2m'-2}e^{i\sqrt{3}\xi_q(m'-1)}\bigg)\\
&& \displaystyle -\bigg(A_2(0,\xi_q)\sum_{n',m'}c_{n'-n,2m'-3}e^{i\sqrt{3}\xi_q(m'-1)}+
A_1(0,\xi_q)\sum_{n',m'}c_{n'-n,2m'-2}e^{i\sqrt{3}\xi_q(m'-1)}\bigg). \label{eq:Delta2}
\end{eqnarray}
\end{normalsize}

Taking Taylor expansion of varibles in terms of $\xi_q$, we can obtain
\begin{eqnarray}
\label{eq:taylor1}\omega(0,\xi_q)&=&\xi_q+o(\xi_q^2),\\
A_1(0,\xi_q)&=&3+i2\sqrt{3}\xi_q-3\xi_q^2+o(\xi_q^2),\\\label{eq:taylor2}
A_2(0,\xi_q)&=&3-\xi_q^2+o(\xi_q^2),\\\label{eq:taylor3}
e^{i\sqrt{3}\xi_q(m'-1)}&=&1+i\sqrt{3}\xi_q(m'-1)-\frac{3(m'-1)^2}{2}\xi_q^2+o(\xi_q^2).\label{eq:taylor4}
\end{eqnarray}

Then substitute Eq. (\ref{eq:taylor1})-(\ref{eq:taylor4}) into Eq. (\ref{eq:Delta2}),  and the residual function can be written as
\begin{eqnarray}
\Delta(0,\xi_q)=D_0+D_1\xi_q+D_2\xi_q^2+o(\xi_q^2).
\end{eqnarray}
Here, $D_\alpha$ $(\alpha=0,1,2)$ are polynomials composed of coefficients $b_{i,j}$ and $c_{i,j}$.

We require $\Delta(0,\xi_q) = o(\xi_q^2)$ to force reflection suppression for long normal incident waves. This yields three linear equations
\begin{eqnarray}
D_\alpha=0, ~~(\alpha=0,1,2). \label{eq:Di}
\end{eqnarray}

Furthermore, in order to suppress reflections of incident waves with broad band wave numbers and big incident angles,  several typical wave vectors $(\xi_p^{\ast},\xi_q^{\ast})$  are taken where the imaginary part and the real part of the residual function  equal zero. That is
\begin{numcases}
{ } {\rm Re} \{\Delta(\xi_p^{\ast},\xi_q^{\ast})\}=0,\label{eq:Re}\\
{\rm Im} \{\Delta(\xi_p^{\ast},\xi_q^{\ast})\}=0. \label{eq:Im}
\end{numcases}

It means that incident wave propagates transparently across the boundary at such a wave vector $(\xi_p^{\ast},\xi_q^{\ast})$. The choice of the wave vector $(\xi_p^{\ast},\xi_q^{\ast})$ is important. It requires a few tries to get an appropriate combination for practical applications.

As the process of solving boundary coefficients is lengthy,  we only put the detailed algebra of the MBC1 as an illustration in text. The details of  the derivation of other matching boundary conditions are given in Appendix A \ref{appen:DetailMBC}.

First, we only use two atoms inside to construct the MBC1, which is the simplest and lowest order of form for the honeycomb lattice as shown in Fig. \ref{f:boundary-4} (a). The MBC1 has three undetermined coefficients, namely $b_{1,1}, c_{0,0}$ and $c_{1,1}$. Thus we just need to require $\Delta(0,\xi_q) = o(\xi_q^2)$, which leads
to three linear equations
\begin{eqnarray} \label{eq:Eqs_LongLimit}
\begin{cases}
D_0=c_{0,0}+2c_{1,1}=0,\\
D_1=3+6b_{1,1}-2\sqrt{3}c_{0,0}-6\sqrt{3}c_{1,1}=0,\\
D_2=2\sqrt{3}+6\sqrt{3}b_{1,1}-3c_{0,0}-11c_{1,1}=0.
\end{cases}
\end{eqnarray}
Solving these equations, we obtain the coefficients
\begin{eqnarray}
b_{1,1}=0.5,~~c_{0,0}=-3.4641,~~c_{1,1}=1.7321.
\end{eqnarray}
Thus the MBC1 is expressed by
\begin{eqnarray}
\dot{v}_{n,1}+0.5(\dot{w}_{n-1,2}+\dot{w}_{n+1,2})=
-3.4641v_{n,1}+1.7321(w_{n-1,2}+w_{n-1,2}).
\end{eqnarray}

Then we choose two layers of atoms to construct the MBC2 as shown in Fig. \ref{f:boundary-4} (b), which contains five undetermined coefficients. In addition to three equations  Eq. (\ref{eq:Eqs_LongLimit}) at long wave limit, we need two additional equations, that is Eq. (\ref{eq:Re}) and Eq. (\ref{eq:Im}) at the wave vector $(\xi_p^{\ast},\xi_q^{\ast})=(0,0.5)$. Solving these five linear coefficient equations, we obtain the MBC2
\begin{eqnarray}\nonumber
&&~~~\dot{v}_{n,1}+0.1671(\dot{w}_{n-1,2}+\dot{w}_{n+1,2})
+0.014(\dot{v}_{n-1,2}+\dot{v}_{n+1,2}) \\
&&=-2.6630v_{n,1}+1.4075(w_{n-1,2}+w_{n-1,2})-0.0760(v_{n-1,2}+v_{n-1,2}).
\end{eqnarray}

We find that the MBC3 has no significant advantages compared with the MBC1 and MBC2. So we do not give the MBC3 in the main body, and put the derivation and analysis in Appendix B \ref{appen:MBC3}.

To further improve the effectiveness of the matching boundary condition, we use more atoms to construct the high order MBC4. The atomic configuration is displayed in Fig. \ref{f:boundary-4}(c). %To make the form of the MBC4 more compact, we let the coefficients of $v_{n,3}$ and $w_{n,3}$ equal to twice the coefficients of $v_{n-2,3}$ and $w_{n-2,3}$ respectively.
Considering the symmetry of the MBC4, it has nine undetermined coefficients. Thus we need to select three wave vectors, namely $(0, 0.5)$, $(0, 1)$ and $(0, 2.4)$  to yields six linear equations. Among these wave vectors, the group velocity has a maximum at $(0, 2.4)$ for the optical branch. Solving three equations at long wave limit Eq. (\ref{eq:Eqs_LongLimit}) and six equations for selected wave vectors simultaneously, we obtain the MBC4
\begin{eqnarray}\nonumber &&~~~\dot{v}_{n,1}+1.6261(\dot{w}_{n-1,2}+\dot{w}_{n+1,2})+0.5183(\dot{v}_{n-1,2}+\dot{v}_{n+1,2}) \\ \nonumber
&&~~~+1.0146(\dot{w}_{n-2,3}+2\dot{w}_{n,3}+\dot{w}_{n+2,3})+0.1034(\dot{v}_{n-2,3}+2\dot{v}_{n,3}+\dot{v}_{n+2,3}) \\ \nonumber
&&=-4.3763 v_{n,1}+0.7035(w_{n-1,2}+w_{n+1,2})+0.2887(v_{n-1,2}+v_{n+1,2}) \\
&&~~~+0.0729(w_{n-2,3}+2w_{n,3}+w_{n+2,3})+0.5250(v_{n-2,3}+2v_{n,3}+v_{n+2,3}). \label{eq:mbc4}
\end{eqnarray}

\begin{figure}
  \centering
  \includegraphics[width=8cm]{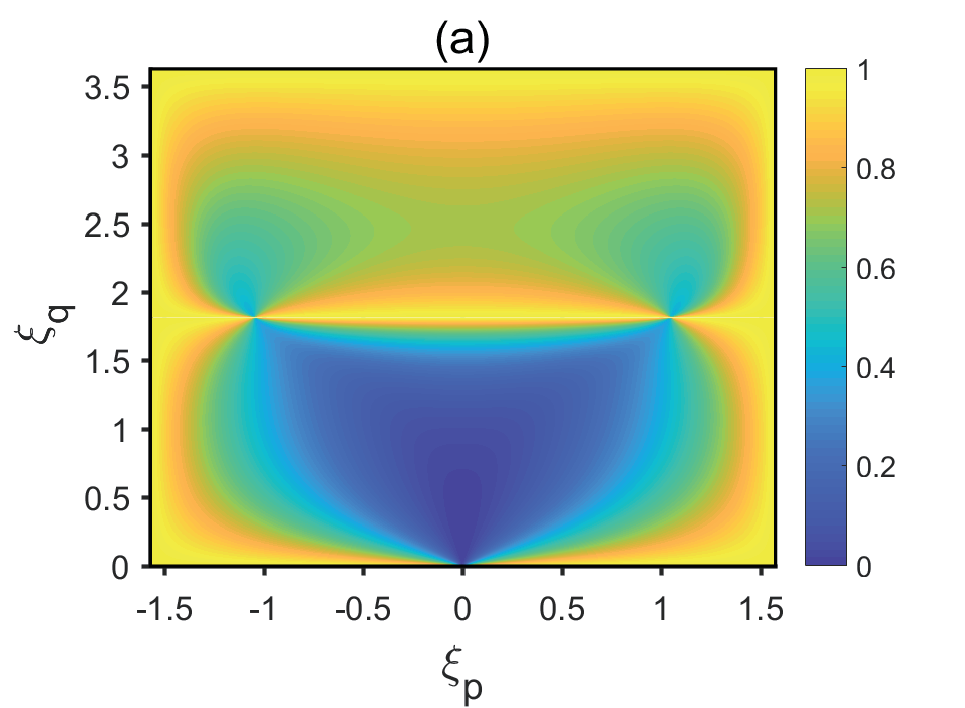}%
  \includegraphics[width=8cm]{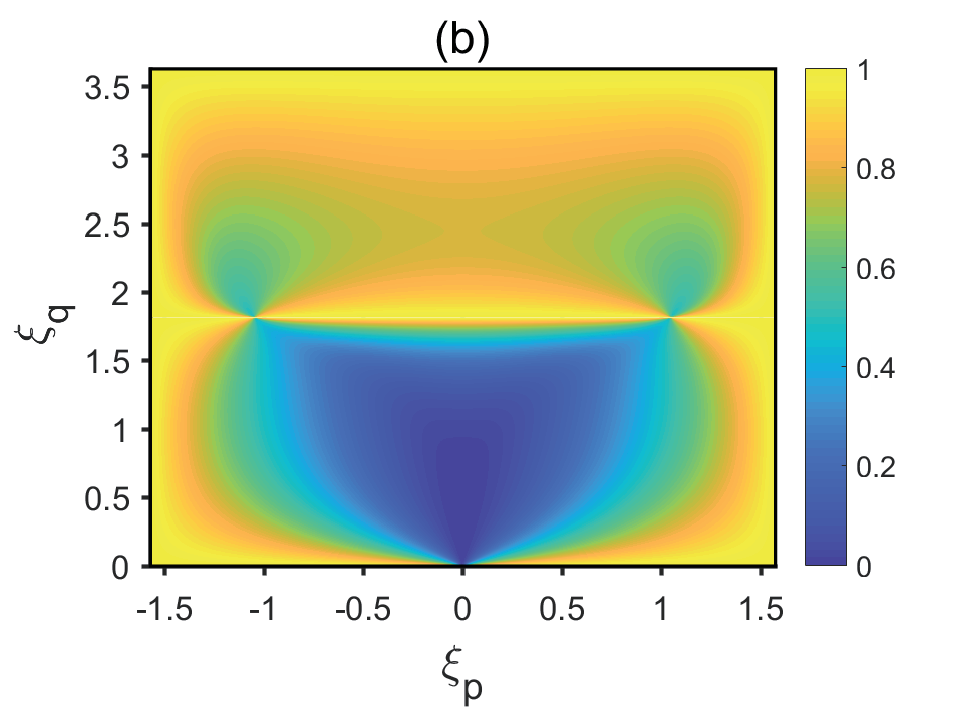}\\
  \includegraphics[width=8cm]{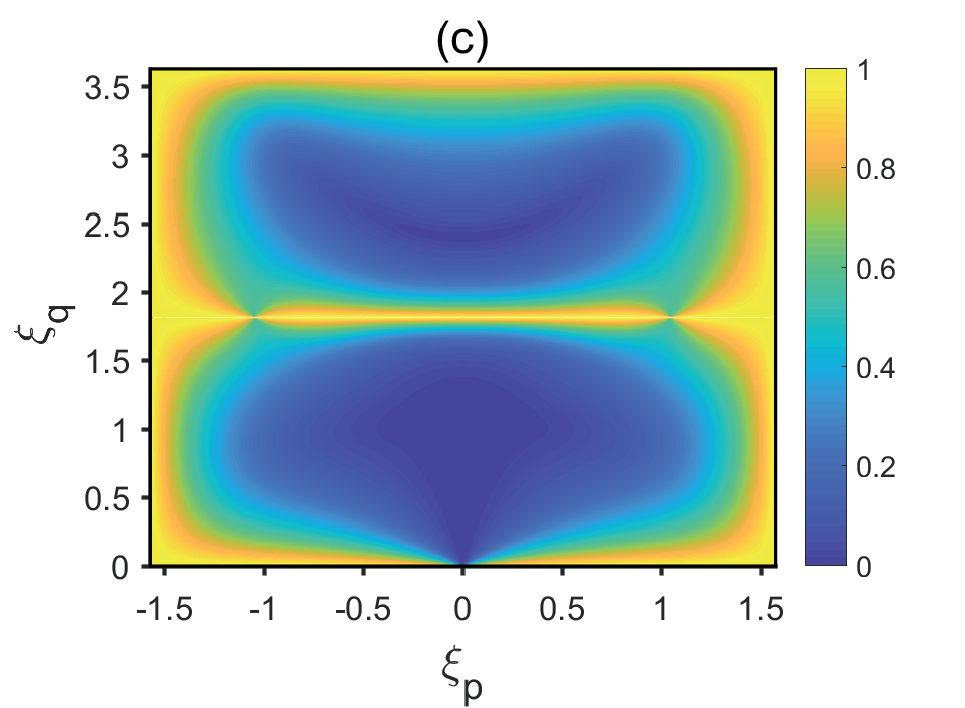}
  \includegraphics[width=8cm]{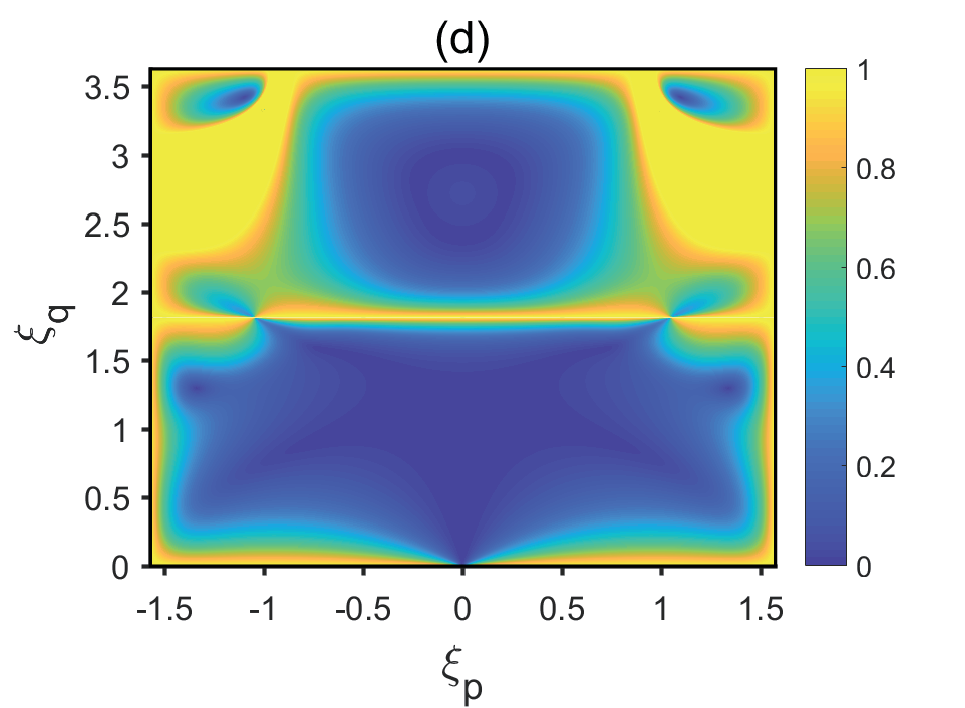}
  \caption{The modulus of reflection coefficients for matching boundary conditions: (a) MBC1; (b) MBC2; (c) MBC4; (d) MBC5.}\label{f:reflection}
\end{figure}

\begin{figure}
  \centering
  \includegraphics[width=8cm]{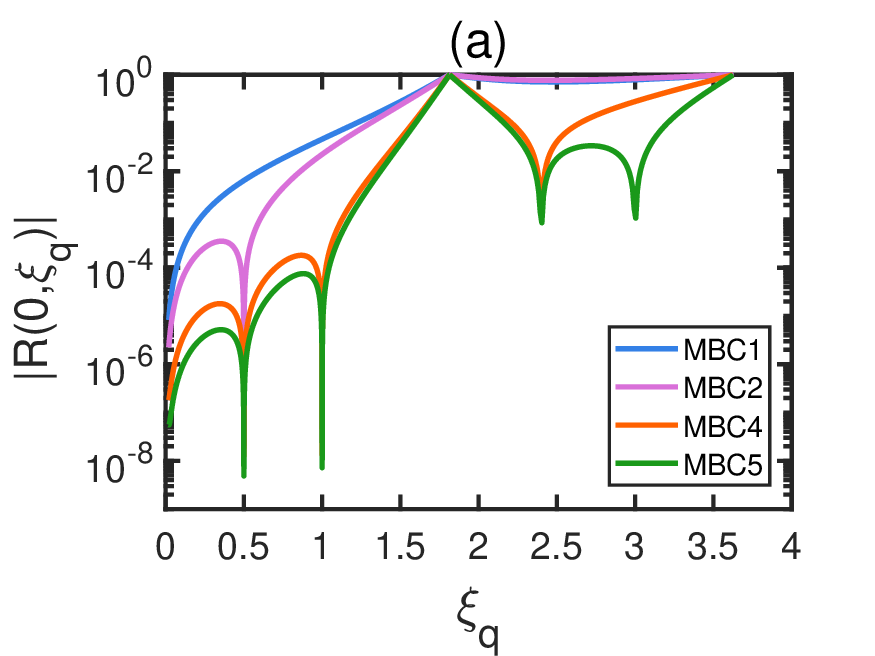}%
  \includegraphics[width=8cm]{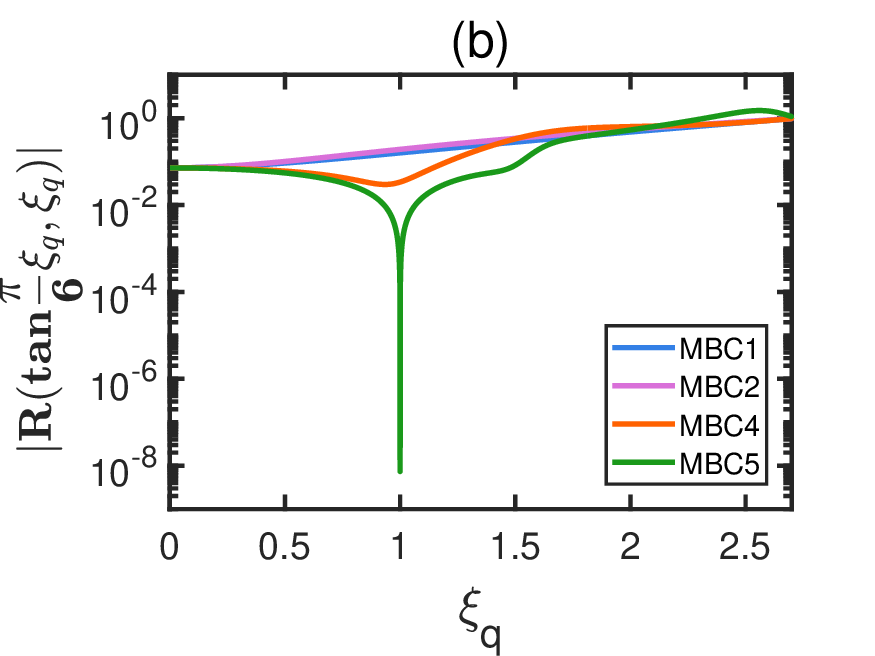}
  \caption{The modulus of reflection coefficients along two special directions: (a) normal direction; (b) oblique direction with an incident angle $\frac{\pi}{6}$.}\label{f:reflection2}
\end{figure}

At last, we add one more layer of atoms to construct the higher order MBC5 as shown in Fig. \ref{f:boundary-4} (d). It contains thirteen undetermined coefficients. Based on the selected wave vectors of the MBC4, we add two additional wave vectors $(\tan\frac{\pi}{6},1)$ and $(0,3)$. The former one is used to expand the oblique incident angles, and the latter one is used to eliminate extremely short waves for the optical branch. After solving corresponding equations, we obtain the MBC5
\begin{eqnarray}\nonumber &&~~~\dot{v}_{n,1}+4.1088(\dot{w}_{n-1,2}+\dot{w}_{n+1,2})+5.7272(\dot{v}_{n-1,2}+\dot{v}_{n+1,2}) \\ \nonumber
&&~~~+2.8636(\dot{w}_{n-2,3}+2\dot{w}_{n,3}+\dot{w}_{n+2,3})+2.0544(\dot{v}_{n-2,3}+2\dot{v}_{n,3}+\dot{v}_{n+2,3}) \\ \nonumber
&&~~~1.0616(\dot{w}_{n-1,4}+\dot{w}_{n+1,4})-0.5616(\dot{w}_{n-3,4}+\dot{w}_{n+3,4})
 \\ \nonumber
&&=-6.4564v_{n,1}-3.2987(w_{n-1,2}+w_{n+1,2})+3.2857(v_{n-1,2}+c_{n+1,2}) \\ \nonumber
&&~~~-1.6428(w_{n-2,3}+2w_{n,3}+w_{n+2,3})+1.6493(v_{n-2,3}+2v_{n,3}+v_{n+2,3}) \\
&&~~~3.7441(w_{n-1,4}+w_{n+1,4})-0.5159(w_{n-3,4}+w_{n+3,4}).
\label{eq:mbc}
\end{eqnarray}

\subsection{2.3 Reflection analysis of matching boundary conditions}

 The effectiveness of an artificial boundary condition may
be checked by reflections. Considering a wave in the form of $v_{n,m} \sim A_1 e^{i(\omega t + \xi_p n + \sqrt{3}\xi_q m)}$ and $w_{n,m} \sim A_2 e^{i(\omega t + \xi_p n + \sqrt{3}\xi_q m)}$ reaching the bottom boundary, reflections will occur due to the artificial boundary conditions. With reflections, the full wave field expression is
\begin{eqnarray}
    \label{eq:FullWave1}
    v^{\text{Full}}_{n,m} &=& A_1 e^{i(\omega t + \xi_p n + \sqrt{3}\xi_q m)} + R A_1 e^{i(\omega t + \xi_p n - \sqrt{3}\xi_q m)} \\
    \label{eq:FullWave2}
    w^{\text{Full}}_{n,m} &=& A_2 e^{i(\omega t + \xi_p n + \sqrt{3}\xi_q m)} + R A_2 e^{i(\omega t + \xi_p n - \sqrt{3}\xi_q m)}
\end{eqnarray}
Substituting Eq. (\ref{eq:FullWave1}) and  Eq. (\ref{eq:FullWave2}) into the matching boundary conditions, we compute the reflection coefficient
\begin{equation}\label{eR}
\hspace{-1.5cm}
R(\xi_p,\xi_q)=-\frac{\Delta(\xi_p,\xi_q)}{\Delta(\xi_p,-\xi_q)}.
\end{equation}
The modulus $|R(\xi_p,\xi_q)|$ directly reflects the performance of matching boundary conditions on suppressing reflections.

We calculate the reflection coefficient and display the modulus $|R(\xi_p,\xi_q)|$ in Fig. \ref{f:reflection}. Darker color represents small boundary reflections. Symmetry about the axis $\xi_p = 0$ is observed, resulting from the choice of atoms and coefficients. It can be seen that the low-order MBC1 and MBC2 mainly suppress incident waves along the normal direction, especially for long wave limit. As the MBC1 and MBC2 have few atoms to match the dispersion residual function, they can not deal with incident waves with big angles for the acoustic branch, or high frequency for the optical branch. In contrast, the high-order MBC4 and MBC5 can do better. They can suppress the boundary reflections for the acoustic branch and the optical branch simultaneously in a broad band. Comparing Fig. \ref{f:reflection}(c) with Fig. \ref{f:reflection}(d), it can be seen that the MBC5 can treat incident waves with higher frequency and more oblique incident angles than the MBC4, due to the MBC5 matching two more selected wave vectors, one along incident angle $\frac{\pi}{6}$, another one close to the right end of the optical branch interval.

Besides, we depict the reflection coefficients along two special directions $(0, \xi_q)$ and $(\tan(\pi/6)\xi_q, \xi_q)$. In Fig. \ref{f:reflection2}, normal incidences are better treated than oblique ones, as more requirements on the dispersion residual are forced along the normal direction. A low reflection is reached at the selected $\xi^{\ast}_q$, as a consequence of our design for the matching boundary conditions.

Therefore, according to different characteristics of the matching boundary conditions, people can choose an appropriate one for application. The low-order MBC1 and MBC2 are simple and efficient, which are suitable for treating long wave problems, such as moving dislocation, or other large deformation in lattice. The high-order MBC4 and MBC5 have complex forms. They can be applied to incident waves with high frequency for defects in lattice, or heat transformation with broad band phonon modes.

\subsection{2.4 Implementations of matching boundary conditions}

\begin{figure}
  \centering
  \includegraphics[width=14cm]{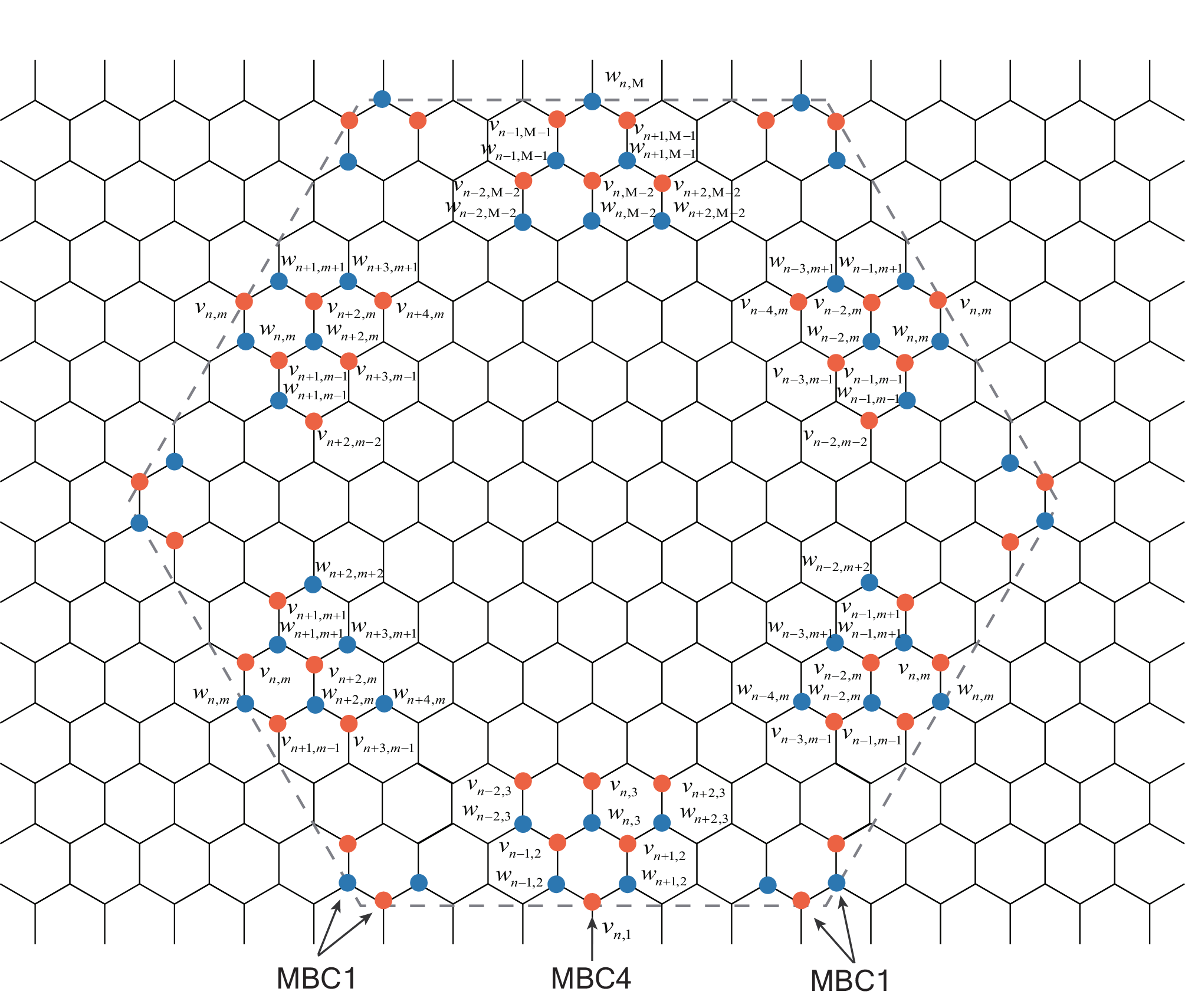}
  \caption{Schematic plot of the hexagon computational domain for the honeycomb lattice with the MBC4 used for six zigzag boundaries and the MBC1 for corner atoms.}\label{f:computing_domain}
\end{figure}

In the following simulations, we extract a hexagon computational domain from a large enough honeycomb lattice with $N=99$ and $M=51$ including 3750 atoms. The schematic plot is shown in Fig. \ref{f:computing_domain}. The matching boundary conditions are applied to the six zigzag boundaries. All corner atoms are taken as the MBC1, represented by the hollow circles in Fig. \ref{f:computing_domain}. When adopting the high-order matching boundary conditions, the boundary conditions of corner atoms need to be treated carefully. This is because the high-order matching boundary conditions involve too many shared atoms in the corner domain, which may cause instability for long time simulations. Previous studies have verified that the low-order matching boundary condition, especially for the MBC1, is stable for long time simulations \cite{S. Ji2014}.  To balance the effectiveness of suppressing reflections and the stability of computations, we combine the high-order matching boundary condition and the MBC1 together in the implementations. The former one is used for ordinary boundary atoms, and the latter one for corner atoms.

Due to different directions of normal lines of the six zigzag boundaries, the matching boundary condition for each boundary has the same coefficients, but different atomic subscripts. In Fig. \ref{f:computing_domain}, we take the MBC4 as an example to illustrate  the form of the matching boundary condition for each zigzag boundary. The MBC4 at the bottom boundary has been given in Eq.(\ref{eq:mbc4}). The MBC4 for other five zigzag boundaries are listed in Appendix C.

\section{3. Numerical tests}

In this section, we first perform atomic simulations for a large enough honeycomb lattice with a harmonic potential Eq. (\ref{eq:harmonic_potential}) to get the reference solution for the propagation of a Gaussian hump. Then we illustrate the effectiveness of matching boundary conditions by comparing the atomic displacements and the kinetic energy with the reference ones. Furthermore, we apply the matching boundary conditions to nonlinear honeycomb lattice with FPU-$\beta$ potential and test their performance of suppressing boundary reflections.

\subsection{3.1 Reflection suppression for harmonic honeycomb lattice}

First, we get the reference solution by simulating a large enough hexagon computational domain and make sure the atomic motions do not reach boundaries.  The velocity verlet scheme are used to integrate the Newton equations with a time step $\Delta t=0.01$.

The initial displacements of $v_{n,m}$ and $w_{n,m}$ are given by $v_{n,m}(0)=u^0(r_{n,m}^{(1)}), w_{n,m}(0)=u^0(r_{n,m}^{(2)})$, with \begin{eqnarray}
u(r)=\Bigg\{
\begin{array}{cc}
A e^{-r^2/50} \cos (r\xi) \cos (\pi r/100), & r \leq 20, \\\\
0, & r>20.
\end{array}
\end{eqnarray}
Here, $r_{n,m}^{(1)}$ and $r_{n,m}^{(2)}$ represent the distance of two classes of atoms numbered by $(n,m)$ from the center of the computational domain respectively, and the amplitude is set $A=0.05$ and the parameter $\xi=0.3$. The velocities of all atoms are set to 0 for all atoms initially.

\begin{figure}
  \centering
  \includegraphics[width=5.8cm]{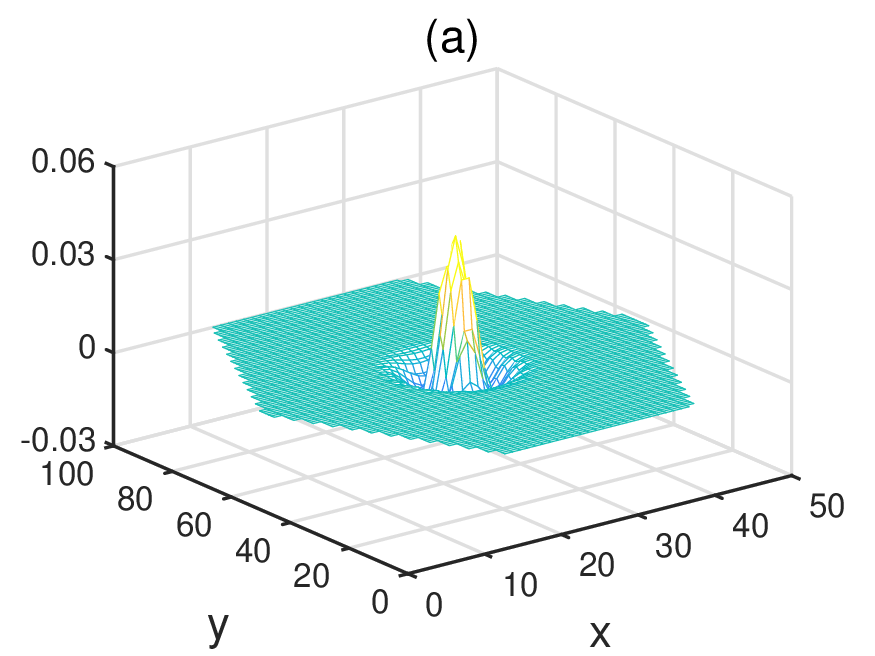}%
  \includegraphics[width=5.8cm]{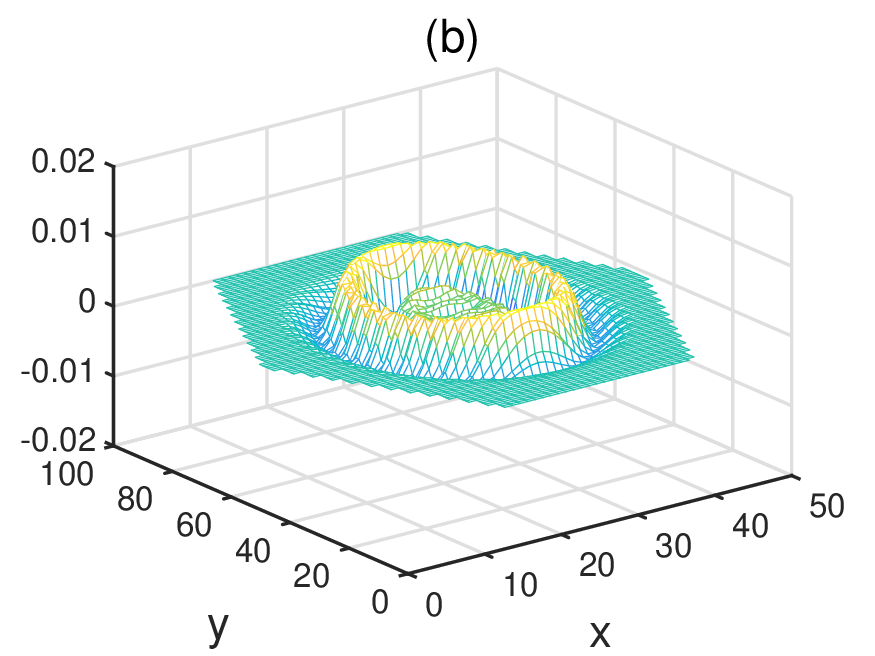}
  \includegraphics[width=5.8cm]{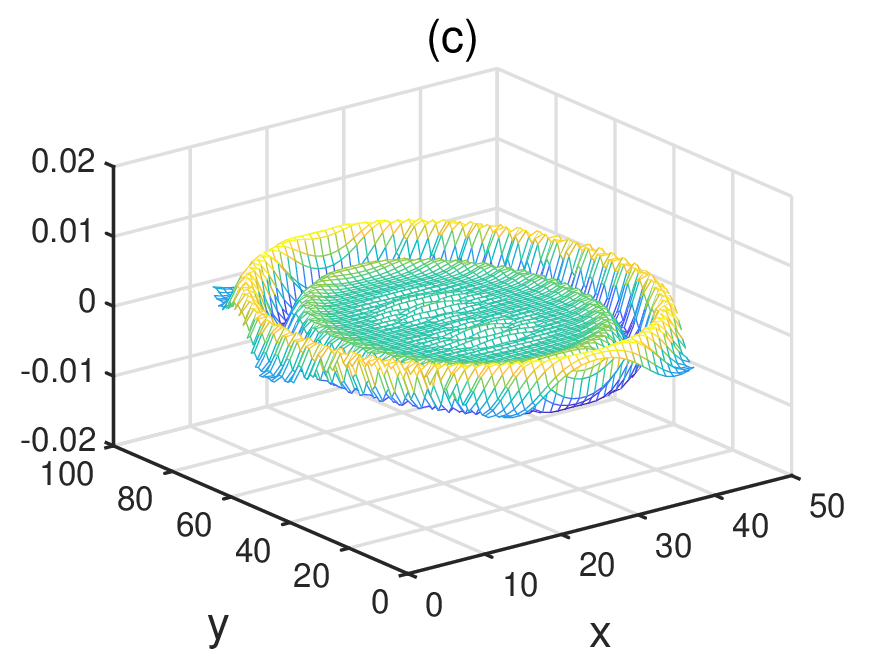}\\
  \includegraphics[width=5.8cm]{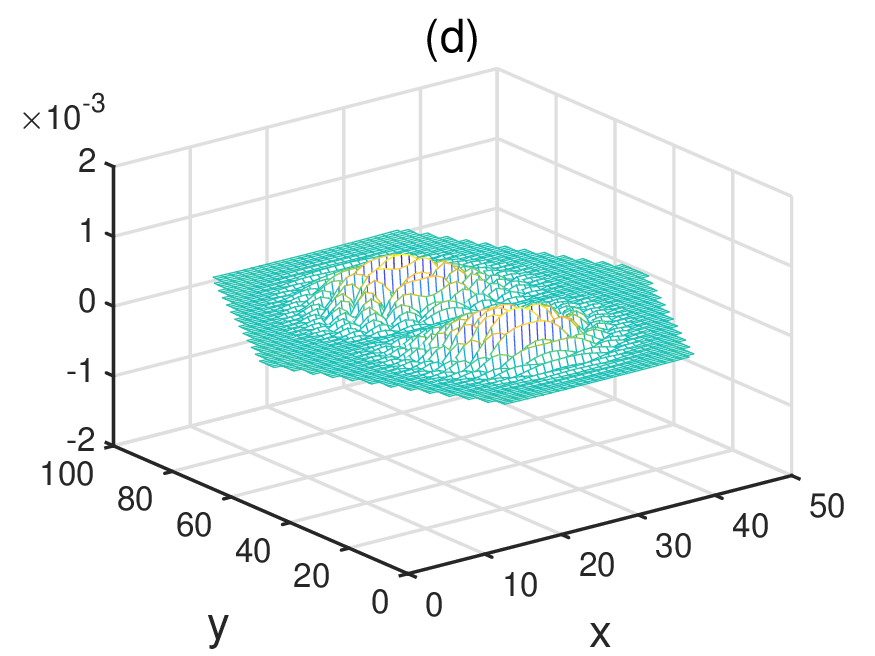}%
  \includegraphics[width=5.8cm]{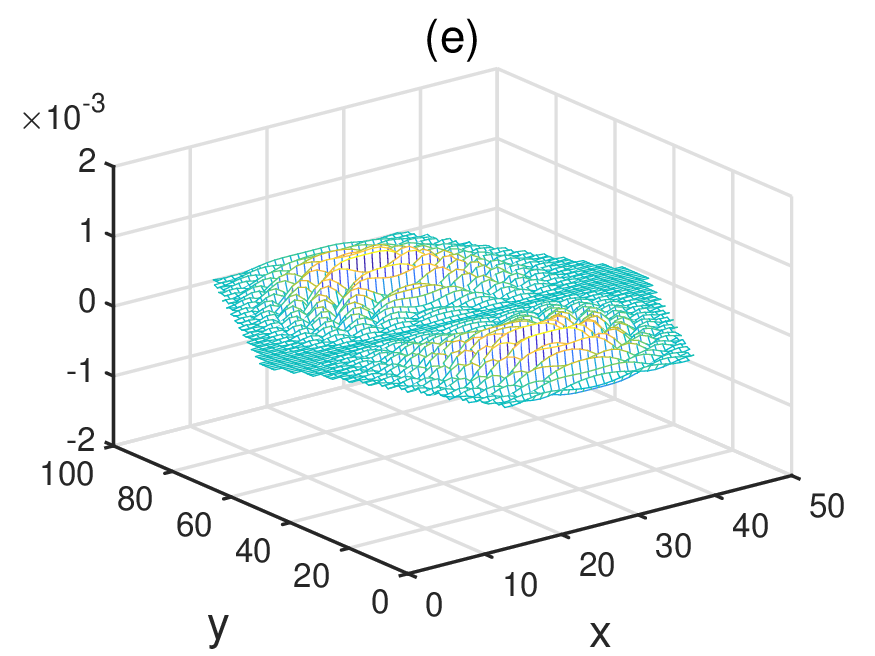}
  \includegraphics[width=5.8cm]{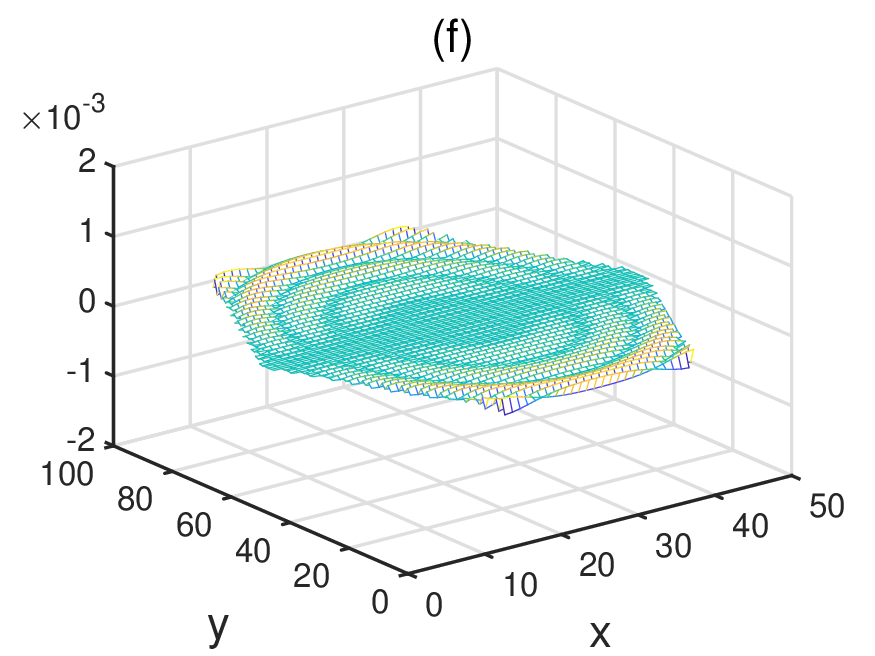}
  \caption{The displacements of the reference solution for the harmonic honeycomb lattice at different times: (a) t=1; (b) t=20; (c) t=40; (d) t=80; (e) t=120; (f) t=300. These figures illustrate a propagation of the Gaussian hump.}\label{f:Guass_harmonic}
\end{figure}
\begin{figure}
  \centering
  \includegraphics[width=6.5cm]{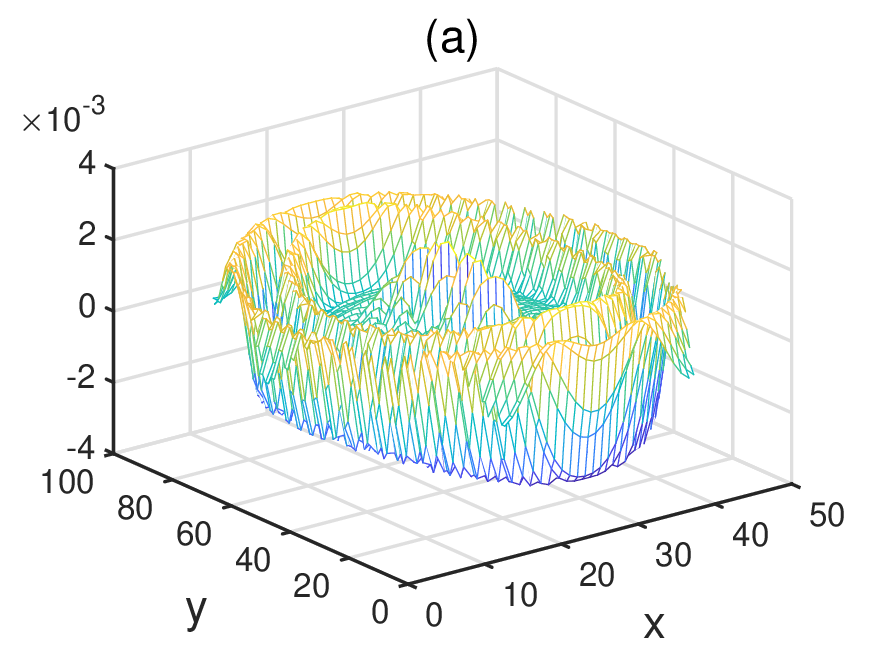}
  \includegraphics[width=6.5cm]{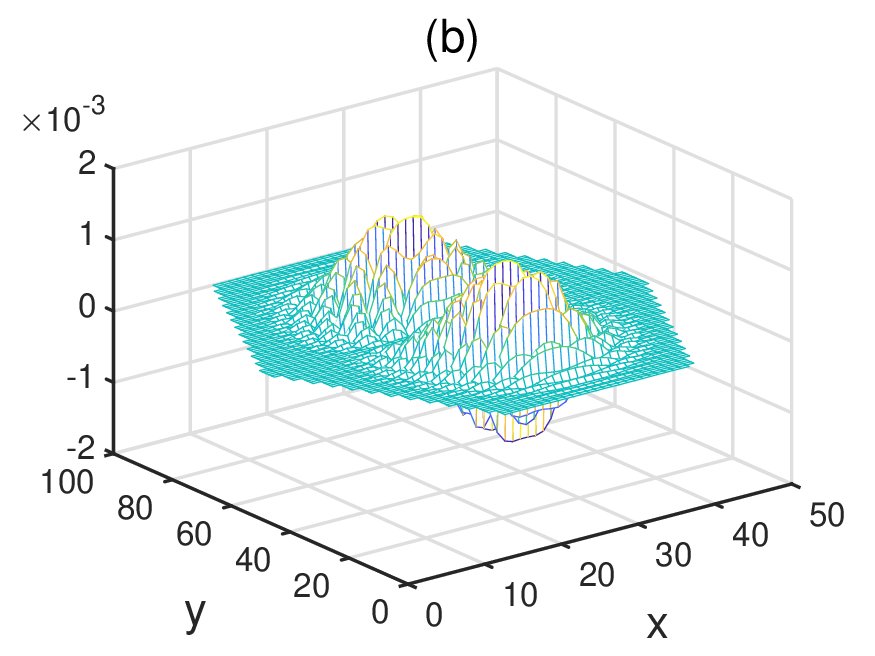}\\
  \includegraphics[width=6.5cm]{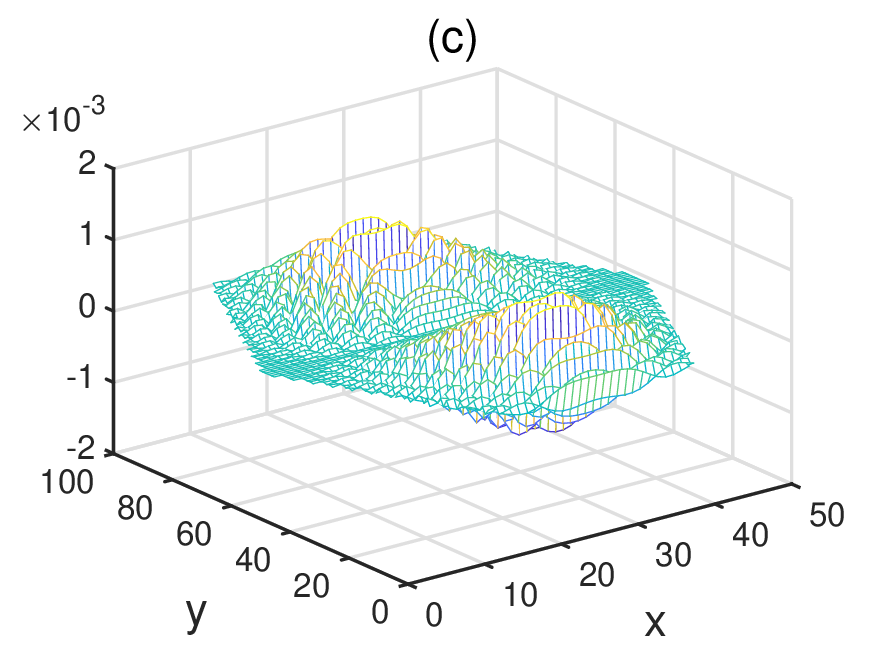}
  \includegraphics[width=6.5cm]{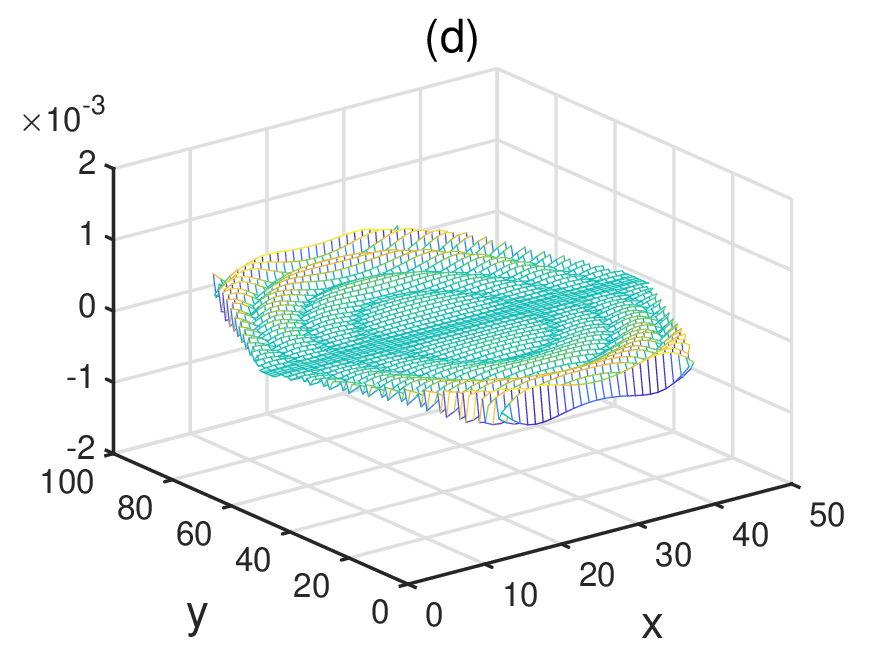}
  \caption{The velocities of the reference solution for the harmonic honeycomb lattice at four moments: (a) t=40; (b) t=80; (c) t=120; (d) t=300.}\label{f:velocity_full}
\end{figure}

We extract a hexagon subdomain with $N=99$ and $M=51$ from a large domain for illustration in the figures. As shown in Fig. \ref{f:Guass_harmonic} (a)-(c), the Gaussian wave packet propagates around gradually. The wave front reaches the boundaries and vanishes gradually after time $t=40$. At time $t=80$ in Fig. \ref{f:Guass_harmonic} (d), the long wave part of the wave packet has passed out of the subdomain, leaving some short waves. These short waves move slowly due to small group velocity as shown in Fig. \ref{f:Guass_harmonic} (e). At time $t=300$ in \ref{f:Guass_harmonic} (f), most of short waves have propagated out of the domain. From the figure, we can also observe that the long waves propagate outward in all directions, while the short waves propagate along the vertical direction.

Simultaneously, we plot the velocities of the honeycomb lattice by mesh grid for the last four moments in Fig. \ref{f:velocity_full}. As we all know, for a monochromatic wave $A e^{i(\omega t + \xi_p n + \sqrt{3}\xi_q m)}$, the amplitude of the velocity is determined by $A\omega$. If a wave has a small amplitude, but owns high frequency, the amplitude of the velocity can be large. Comparing Fig. \ref{f:Guass_harmonic} (c) with Fig. \ref{f:velocity_full} (a), the amplitudes of the center waves are very small, but their velocities are significant in the same order of magnitude as long waves. This indicates that waves at the center are mainly short waves with high frequency. As shown in Fig. \ref{f:velocity_full} (b) and (c), the propagation of short waves can be observed more clearly. At the last moment $t=300$, the displacements of waves are almost invisible in Fig. \ref{f:Guass_harmonic} (f). However, in Fig. \ref{f:velocity_full} (d),  there are small oscillations left in the middle domain. It will take a long time for these short waves passing out of the subdomain completely.

\begin{figure}
  \centering
  \includegraphics[width=6.5cm]{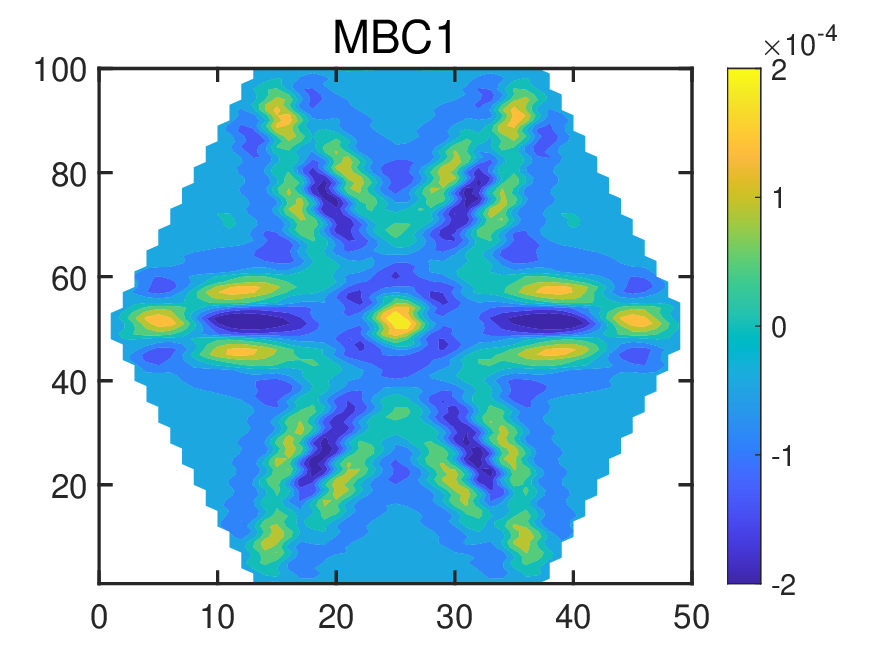}
  \includegraphics[width=6.5cm]{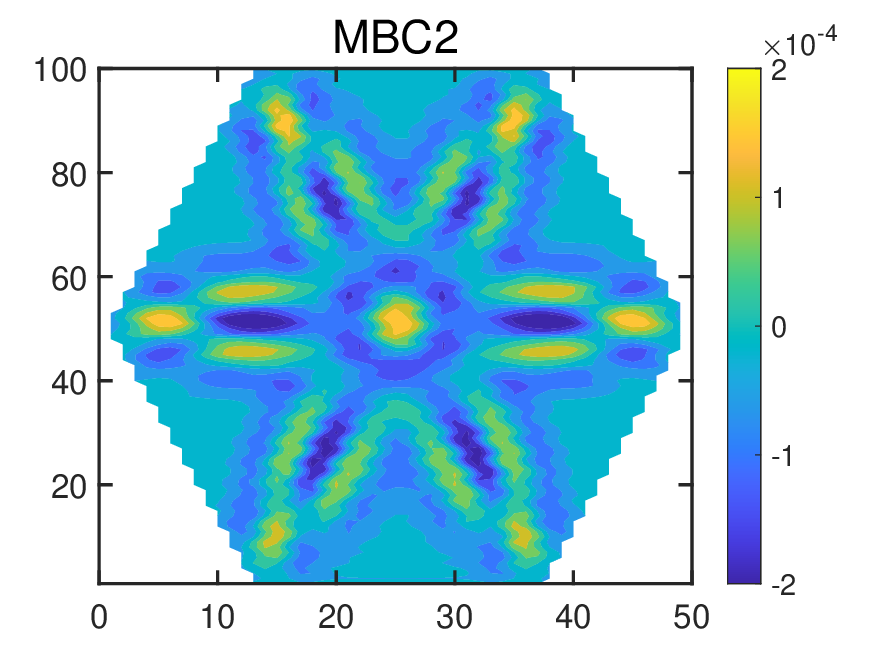}\\
  \includegraphics[width=6.5cm]{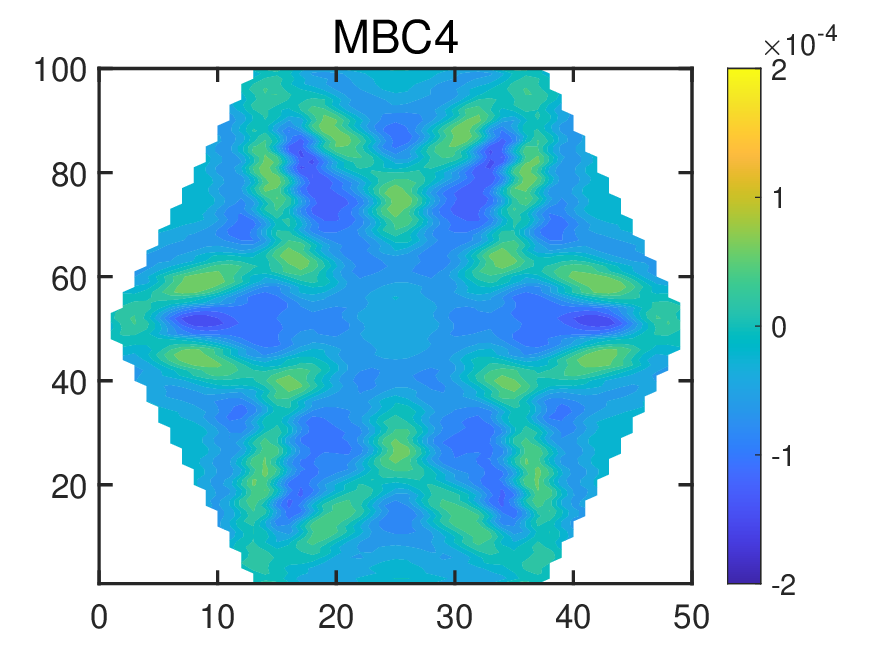}
  \includegraphics[width=6.5cm]{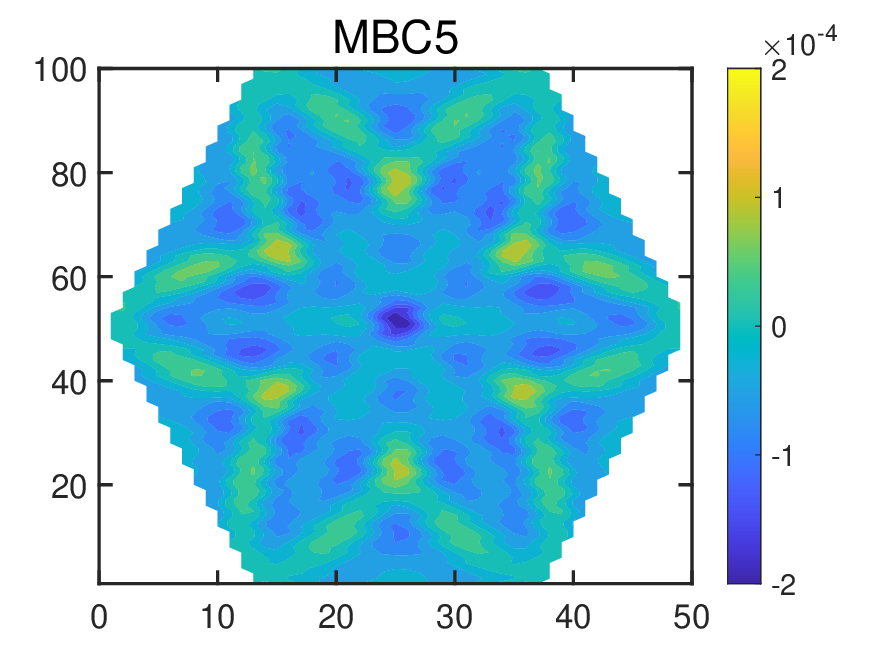}
  \caption{The deviation of displacements $\Delta u(t)$ for the harmonic honeycomb lattice at time $t=100$.}\label{f:t100_harmonic}
\end{figure}
\begin{figure}
  \centering
  \includegraphics[width=6.5cm]{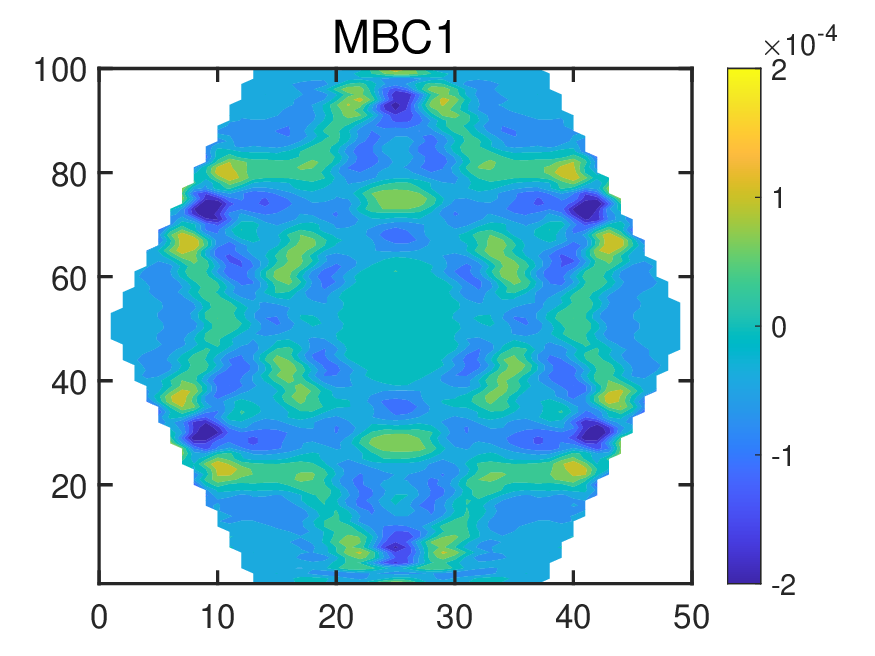}
  \includegraphics[width=6.5cm]{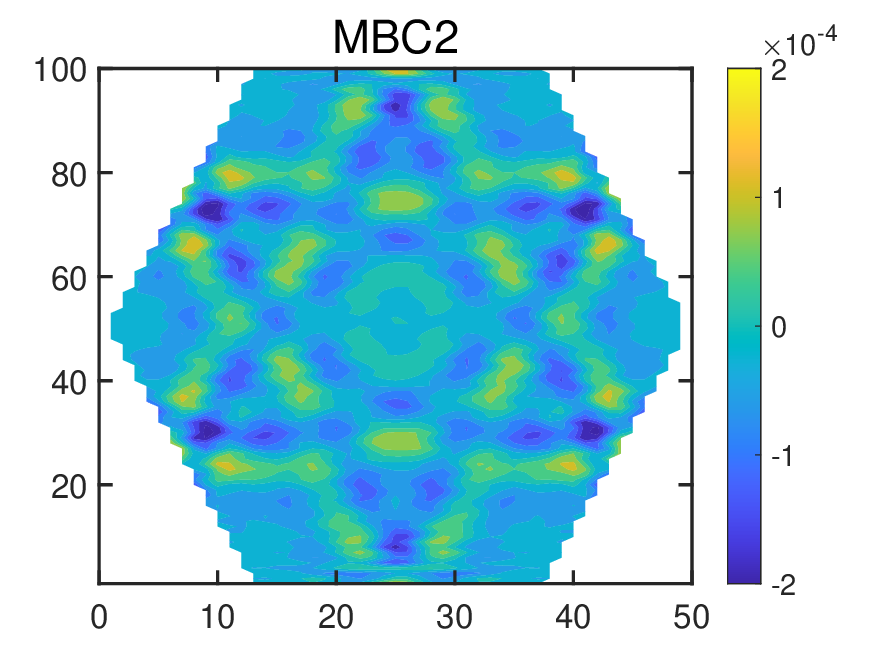}\\
  \includegraphics[width=6.5cm]{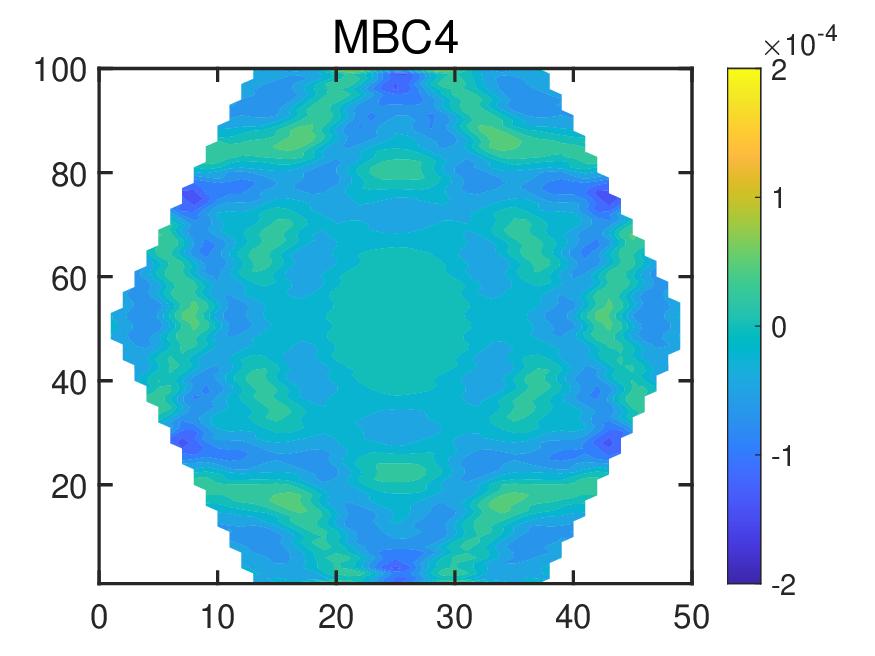}
  \includegraphics[width=6.5cm]{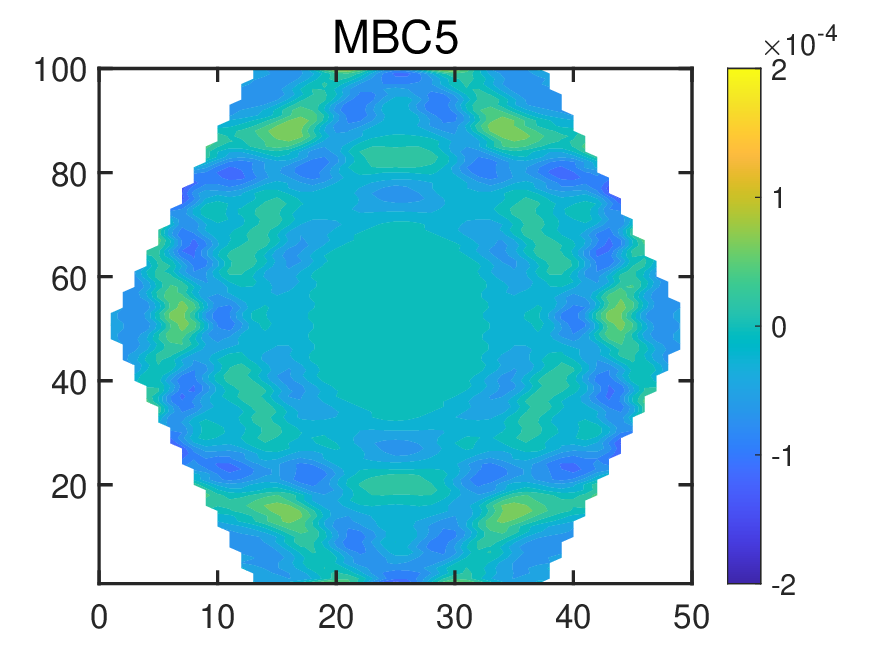}
  \caption{The deviation of displacements $\Delta u(t)$ for the harmonic honeycomb lattice at time $t=120$.}\label{f:t120_harmonic}
\end{figure}
\begin{figure}
  \centering
  \includegraphics[width=6.5cm]{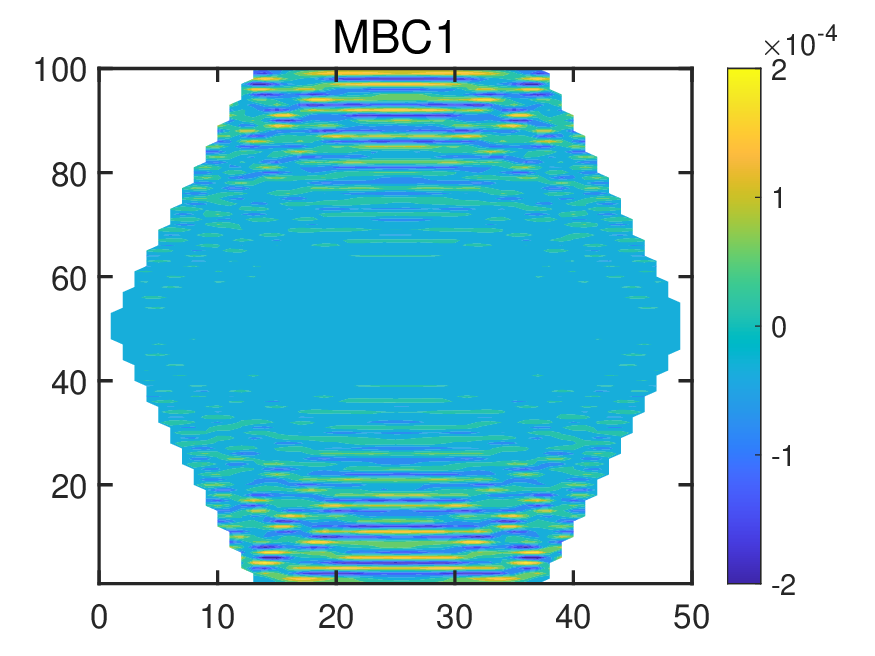}
  \includegraphics[width=6.5cm]{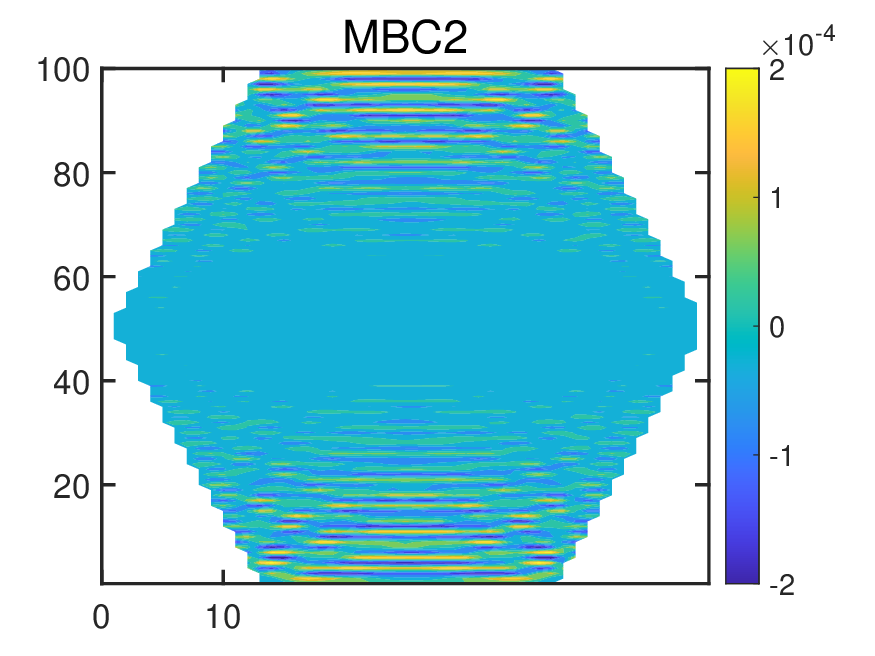}\\
  \includegraphics[width=6.5cm]{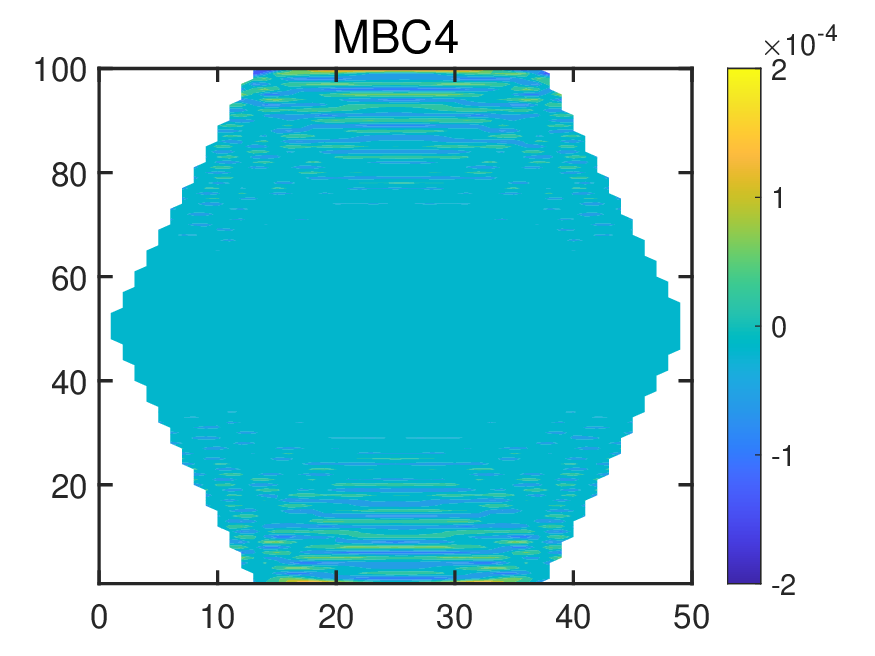}
  \includegraphics[width=6.5cm]{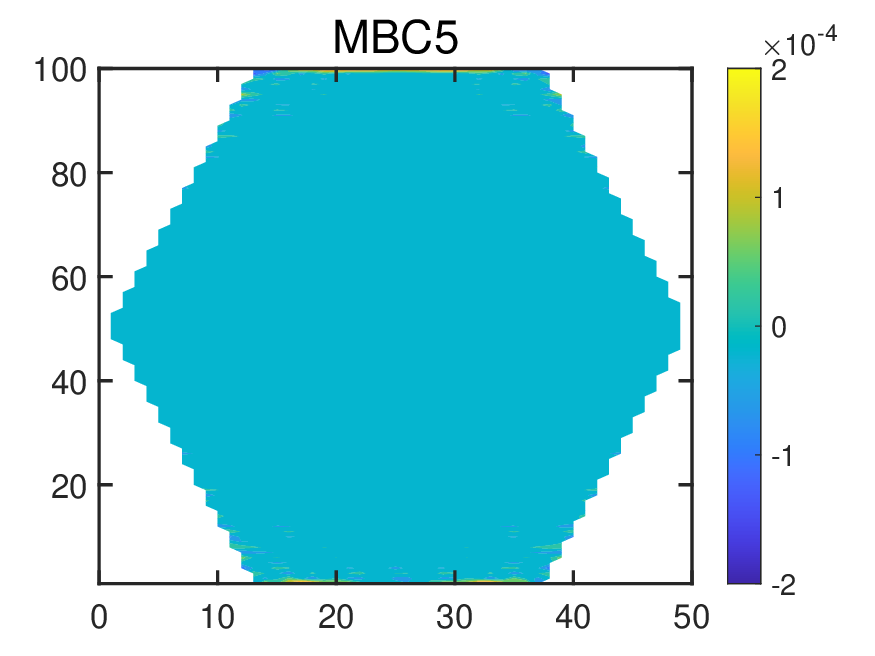}
  \caption{The deviation of displacements $\Delta u(t)$ for the harmonic honeycomb lattice at time $t=200$.}\label{f:t200_harmonic}
\end{figure}

Then we perform atomic simulations for four matching boundary conditions in a hexagon computational domain with $N=99$ and $M=51$ respectively. We set the same initial Gaussian hump in the following simulations. The matching boundary conditions are applied to the six zigzag boundaries with the MBC1 for corner atoms. To clearly show the effectiveness of the matching boundary conditions, we draw the displacement deviation $\Delta u(t)$ of the whole lattice at three different moments in Fig. \ref{f:t100_harmonic}-\ref{f:t200_harmonic}, namely $t=100, t=120, \text{and } t=200$. The deviation $\Delta u(t)$ is defined by
\begin{eqnarray}
\Delta u(t)=\left(
\begin{array}{c}
v_{n,m}(t)-v^{ref}_{n,m}(t)\\
w_{n,m}(t)-w^{ref}_{n,m}(t)
\end{array}
\right).
\end{eqnarray}
Here, $v^{ref}_{n,m}(t)$ and $w^{ref}_{n,m}(t)$ are the reference solution obtained from simulations in a large enough honeycomb lattice.

At time $t=100$, the long waves have propagated out, while the major front of short waves has not yet reached boundaries. At this moment, the boundary reflections are mainly caused by long waves. As shown in Fig. \ref{f:t100_harmonic}, the boundary reflections of the MBC1 and MBC2 are larger than the MBC4 and MBC5. This is because the MBC4 and MBC5 can better deal with oblique incident waves. The patterns are symmetrical, like a six-pointed star, since the long waves spreading out in all directions.

Comparing Fig. \ref{f:t120_harmonic} with Fig. \ref{f:t100_harmonic}, we can see that the color of the pattern at time $t=120$ becomes lighter for each matching boundary condition, illustrating that the reflecting waves vanish gradually. %The color of corner domain is not always dark and the patterns change with time, which means that corner reflections are not the major factor.
According to our observations, the corner reflections are not the major factor, and the boundary reflections mainly come from large angle incident waves. A close inspection, there are some ripples around the bottom and the top boundaries in the images of the  MBC1 and MBC2. This is due to the short waves mainly spreading along the vertical direction. Thus the reflections at the bottom and the top boundaries will be stronger than other four boundaries.

The deviation of atomic displacements at time $t=200$ is displayed in Fig. \ref{f:t200_harmonic}. At this moment, the symmetrical patterns have disappeared, indicating that the reflections caused by long waves have almost been eliminated. There are some fine reflections left in the computational domain. The reflections of short waves  mainly appear around the top boundary and the bottom boundary.  There exist big difference between patterns of these four matching boundary conditions. The ripples of the MBC1 and MBC2 looks the same and obviously distributed along the vertical direction. The ripples of the MBC4 are smaller with the lighter color than the MBC1. In the image of the MBC5, there are no significant ripples, indicating that the MBC5 suppress reflections of short waves effectively.

\begin{figure}
  \centering
  \includegraphics[width=8cm]{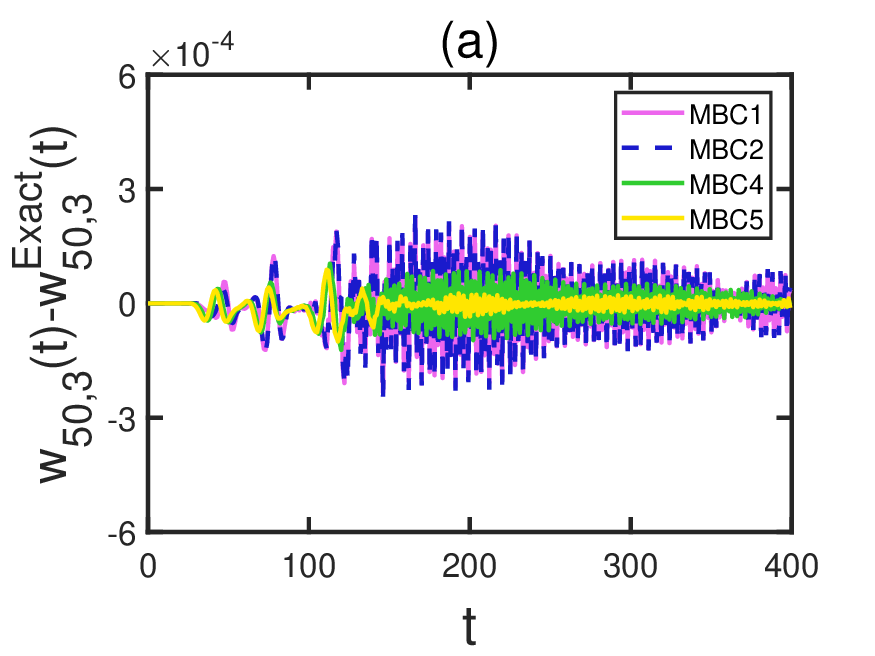}%
  \includegraphics[width=8cm]{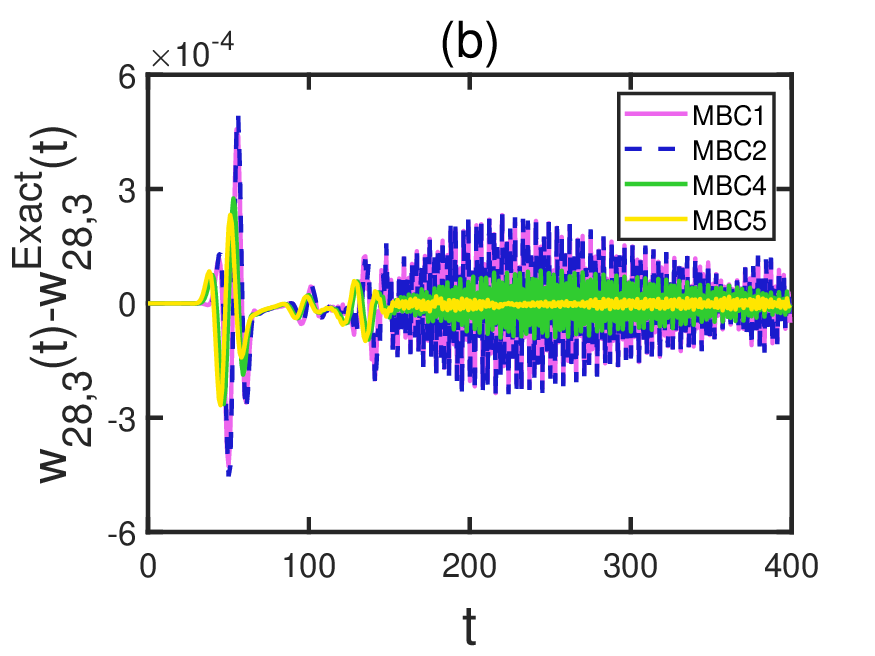}%
  \caption{The deviation of displacements from the reference solution for the harmonic honeycomb lattice: (a) $w_{50,3}(t)$-$w_{50,3}^{ref}(t)$;  (b) $w_{28,3}(t)$-$w_{28,3}^{ref}(t)$.}\label{f:deviation_harmonic}
\end{figure}

\begin{figure}
  \centering
  \includegraphics[width=10cm]{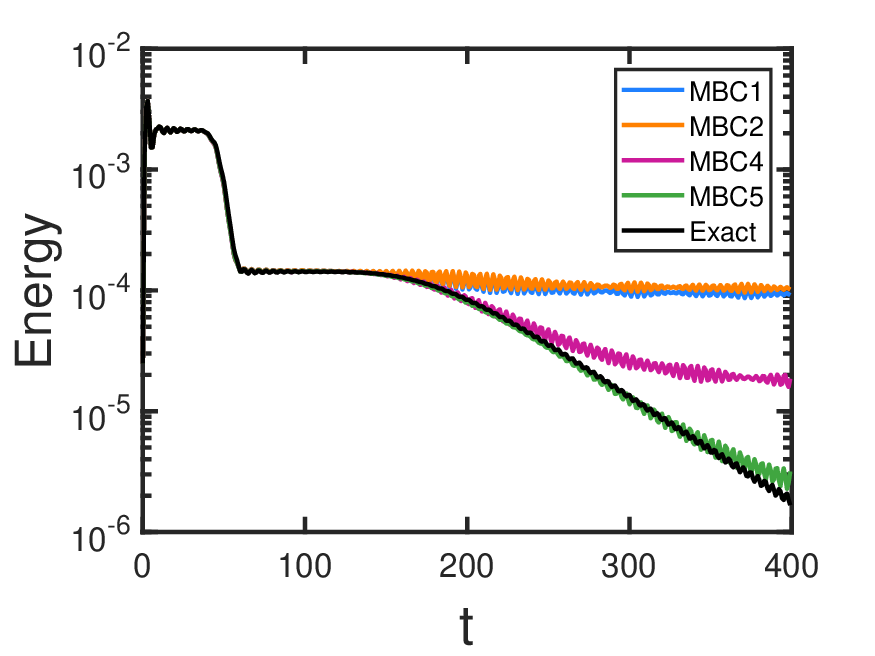}%
  \caption{The kinetic energy of the harmonic honeycomb lattice during the propagation of the Gaussian hump.}\label{f:energy_harmonic}
\end{figure}

To better illustrate the effectiveness of the matching boundary condition along with time, we select two representative atoms and compare their displacements with the  reference solution. The atom $w_{50,3}$ lies at the center of the bottom boundary. Another atom $w_{28,3}$ lies at the lower left corner. The difference $w_{50,3}(t)-w^{ref}_{50,3}(t)$ and $w_{28,3}(t)-w^{ref}_{28,3}(t)$ between the numerical solution and the reference solution are shown in Fig. \ref{f:deviation_harmonic}. The lines of the MBC1 and MBC2 almost coincide, that means the effect of suppressing reflections is similar. The difference among lines mainly appears after time $t=120$ when short waves reach the boundaries. The deviations of the high-order MBC are obviously smaller than the low-order MBC. Especially for the MBC5, the amplitude of the yellow line is quite small,  indicating that the higher order matching boundary condition well suppresses reflections  of short waves. Comparing Fig. \ref{f:deviation_harmonic}(a) with Fig. \ref{f:deviation_harmonic}(b), we observe that the difference in $w_{28,3}(t)$ is bigger than that in $w_{50,3}(t)$ during the long wave passing outward, which is due to the so-called corner reflections. Corner reflection is a subtle issue, and very expensive to resolve cleanly \cite{G. Pang2018 corner}.

At last, we plot the kinetic energy curve of the honeycomb lattice to test the effectiveness and stability of the matching boundary condition. The kinetic energy is defined by
\begin{eqnarray}
W(t)=\sum_{n,m} \frac{1}{2} \Big( \dot{v}^2_{n,m}(t) + \dot{w}^2_{n,m}(t) \Big).
\end{eqnarray}
As shown in Fig. \ref{f:energy_harmonic}, the dark line represents the reference solution. There are two platforms in energy curves. The first platform is corresponding to the process of waves propagating into the lattice. A sudden drop of energy occurs due to the long wave part passing out of the computational domain. Within the time $t=0\sim 140$, all colored curves coincide well with the black one. And then the curves of the MBC1 and MBC2 decline very slowly compared with the reference one after the time $t=140$. As a comparison, the curves of the MBC4 and MBC5 keep decreasing quickly. The absolute value of the slop for the MBC5 is larger than the MBC4. The green curve of the MBC5 is also close to the reference one. This indicates that choosing more atoms to construct matching boundary condition has better effect to eliminate boundary reflections. Furthermore, the energy of these matching boundary conditions do not appear accumulation over time, illustrating the stability of matching boundary conditions numerically.

\subsection{3.2 Reflection suppression for nonlinear honeycomb lattice}

The nonlinear honeycomb lattice is more complicated than harmonic lattice because the nonlinearity will make energy exchange between normal modes. Some short waves with high  frequency will be excited by nonlinear interaction, which makes suppressing boundary reflections difficult. To illustrate the matching boundary conditions for the nonlinear interaction potential, we consider the FPU-$\beta$ potential \cite{S. Lepri2003}
\begin{eqnarray}
U(r)=\sum\frac{k}{2}(r-r_0)^2+\sum\frac{\beta}{4}(r-r_0)^4.
\end{eqnarray}
Here, $r$ is the distance between the nearest atom, $r_0$ is the atomic distance at equilibrium, and $\beta$ is the nonlinear coefficient. Rescale time by $\sqrt{m/k}$ and the displacement by lattice constant $a$. We take the dimensionless $k=1$ and $\beta=1$. Under these scalings, the non-dimensionalised Newton equation with the nearest neighboring interaction takes the form as
\begin{eqnarray}\label{eq:motion_FPU}\nonumber
   \ddot{v}_{n,m}&=&(w_{n-1,m+1}+w_{n+1,m+1}+w_{n,m}-3v_{n,m})\\
   &&+\beta\Big((w_{n-1,m+1}-v_{n,m})^3+(w_{n+1,m+1}-v_{n,m})^3+(w_{n,m}-v_{n,m})^3 \Big), \\\nonumber
   \ddot{w}_{n,m}&=&(v_{n-1,m-1}+v_{n+1,m-1}+v_{n,m}-3w_{n,m})\\
   &&+\beta\Big((v_{n-1,m-1}-w_{n,m})^3+(v_{n+1,m-1}-w_{n,m})^3+(v_{n,m}-w_{n,m})^3\Big).
\end{eqnarray}
\begin{figure}
  \centering
  \includegraphics[width=6.5cm]{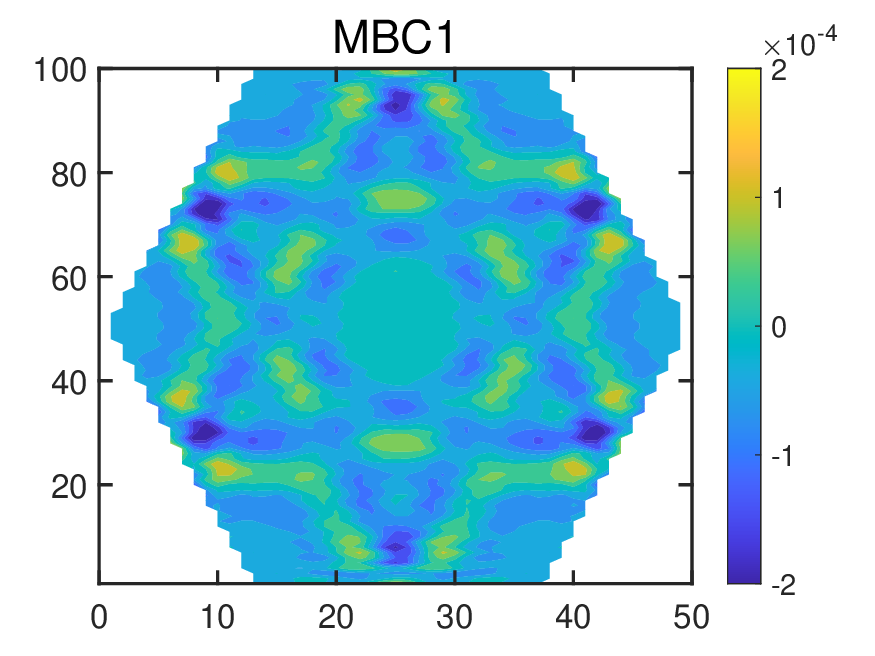}
  \includegraphics[width=6.5cm]{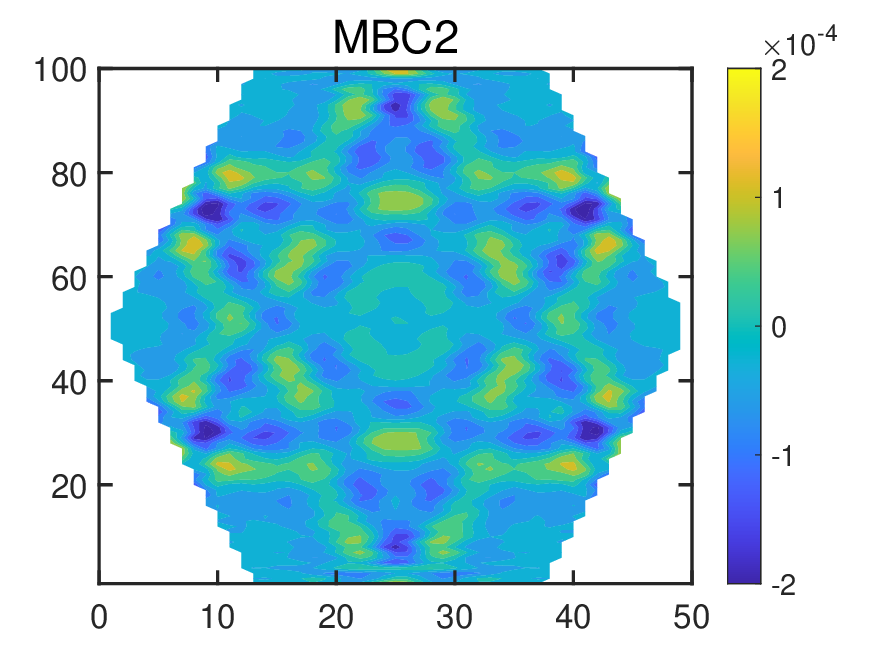}\\
  \includegraphics[width=6.5cm]{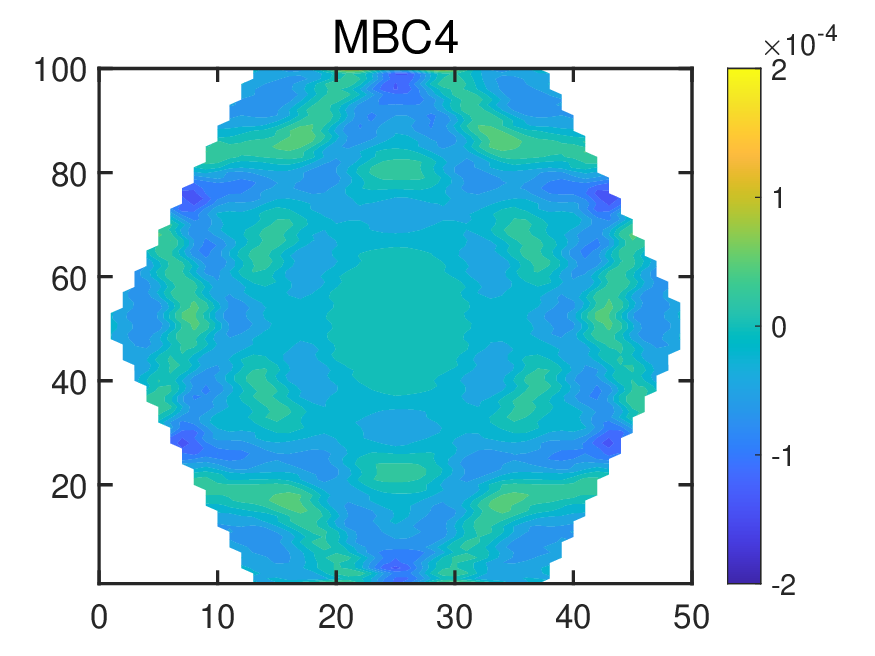}
  \includegraphics[width=6.5cm]{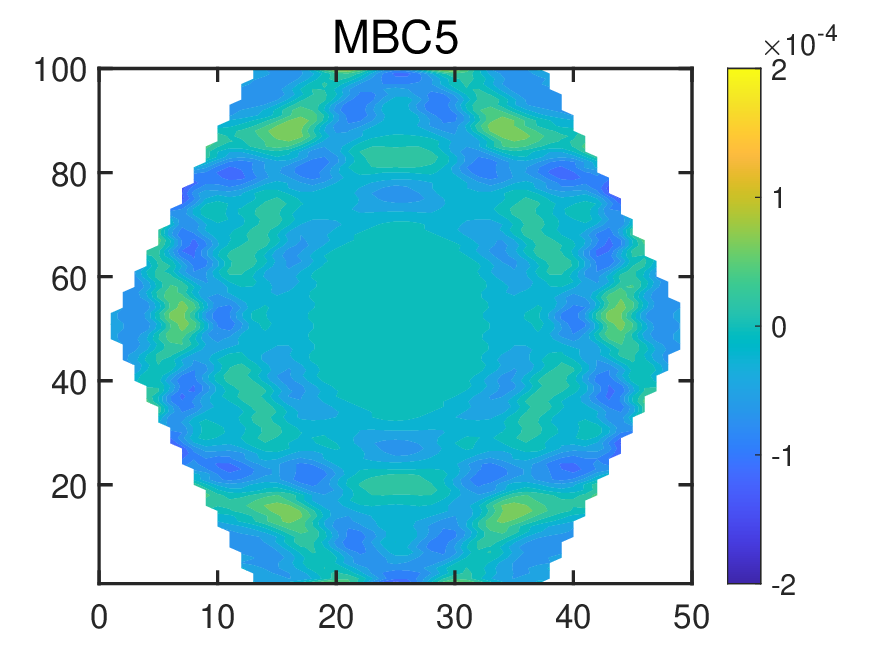}
  \caption{The deviation of displacements for the nonlinear honeycomb lattice at t=120.}\label{f:nonlinear_t120}
\end{figure}
\begin{figure}
  \centering
  \includegraphics[width=6.5cm]{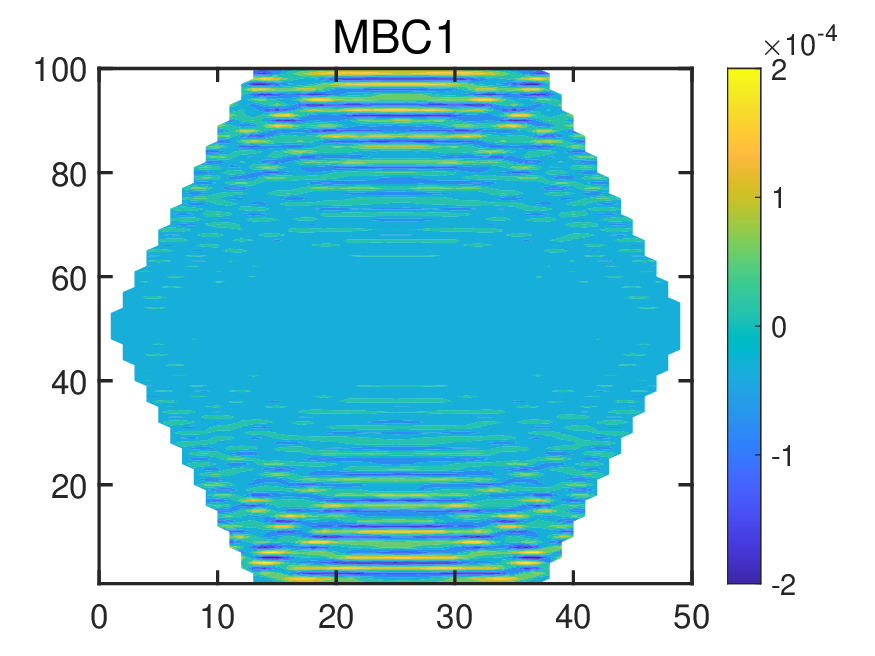}
  \includegraphics[width=6.5cm]{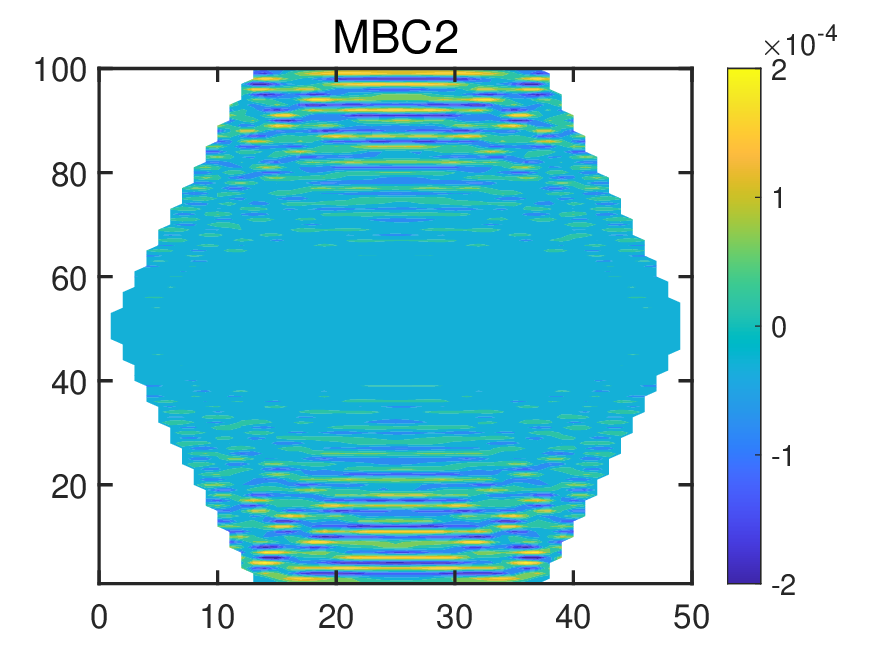}\\
  \includegraphics[width=6.5cm]{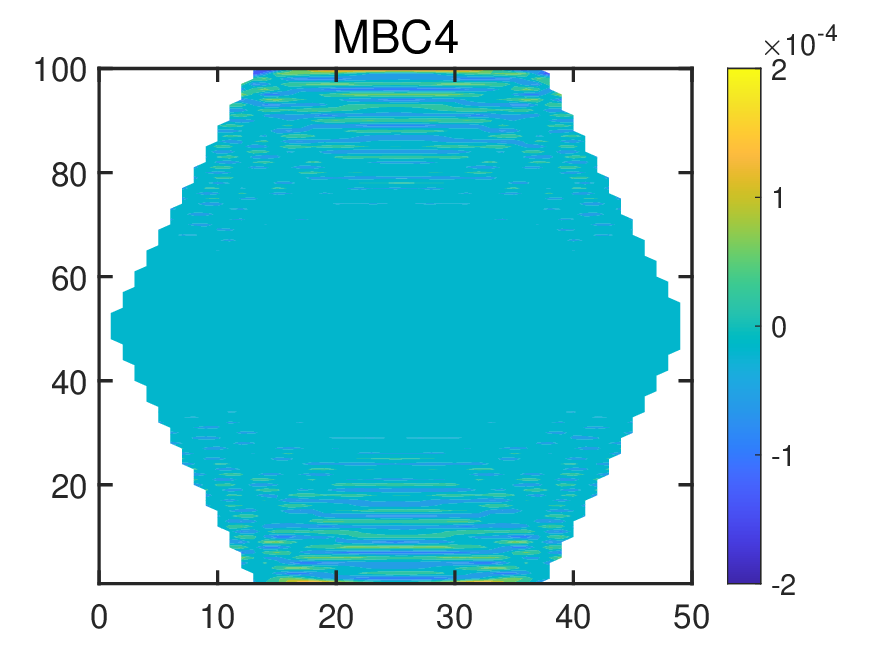}
  \includegraphics[width=6.5cm]{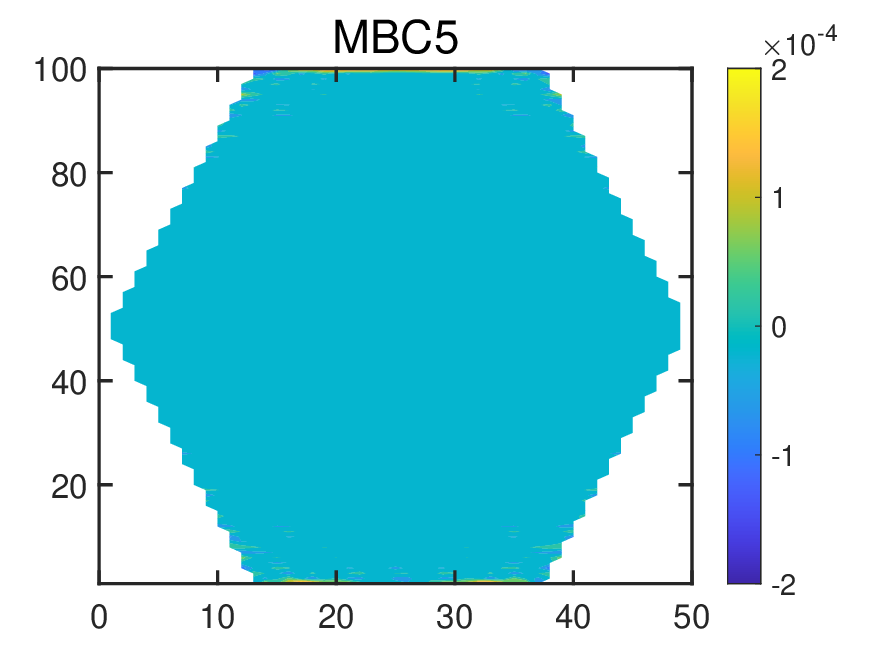}
  \caption{The deviation of displacements for the nonlinear honeycomb lattice at t=200.}\label{f:nonlinear_t200}
\end{figure}

Then we perform a series of atomic simulations for the nonlinear honeycomb lattice with matching boundary conditions. The computational domain and the initial Gaussian hump are taken the same as the harmonic case. The four matching boundary conditions are also applied to the six zigzag boundaries respectively with the MBC1 for the corner atoms. To test applicability and effectiveness, we compare numerical results with the reference solution, which is obtained by simulating a large enough honeycomb lattice with FPU-$\beta$ potential. The deviations of displacement $\Delta u(t)$ at time $t=120$ and $t=200$ are shown in Fig. \ref{f:nonlinear_t120} and Fig. \ref{f:nonlinear_t200}. We can see that the patterns of the nonlinear honeycomb lattice are almost the same as the harmonic cases. It demonstrates that the matching boundary conditions are also effective for the nonlinear honeycomb lattice.

\begin{figure}
  \centering
  \includegraphics[width=8.5cm]{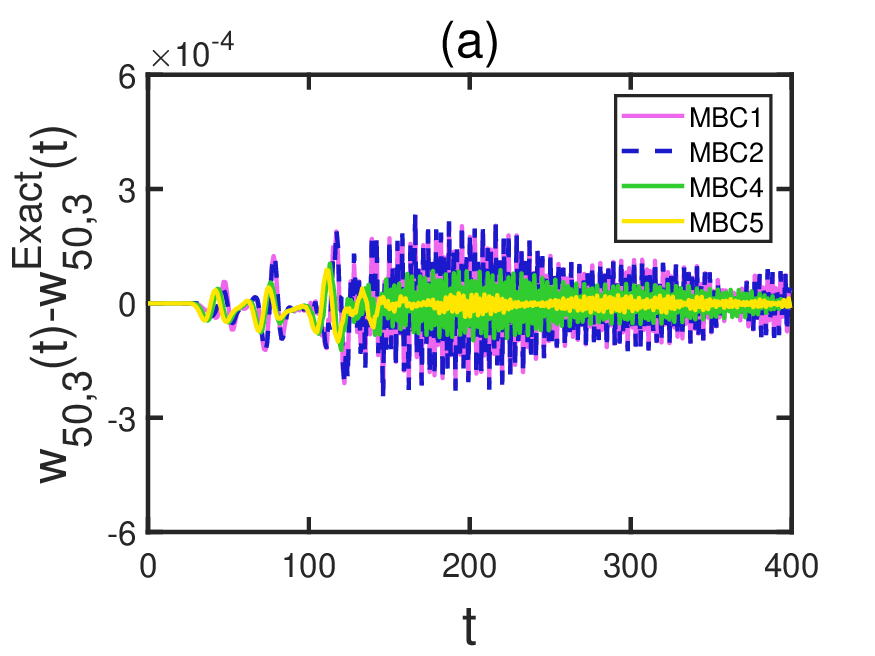}%
  \includegraphics[width=8.5cm]{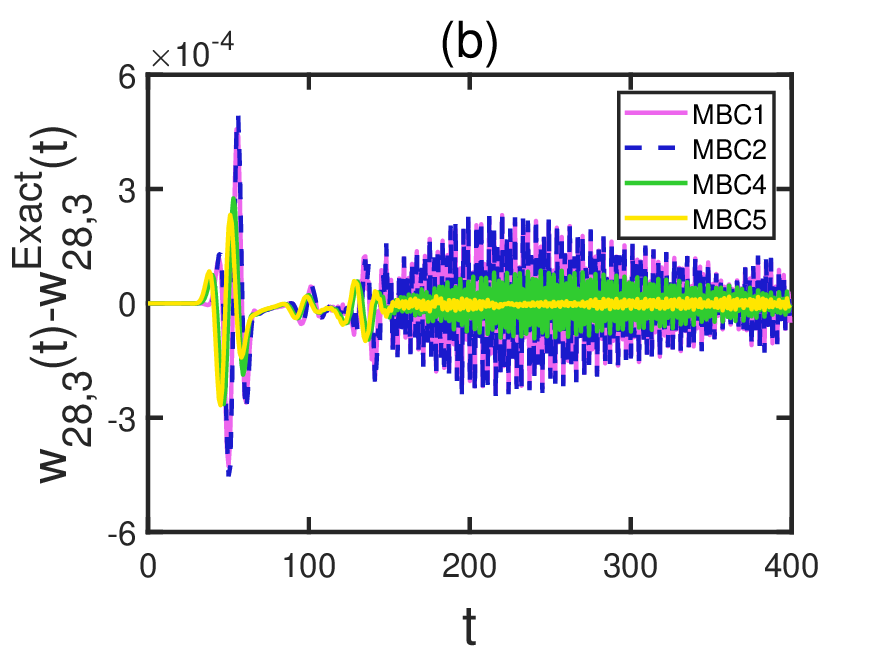}%
  \caption{The deviation of displacements from the reference solution for the nonlinear honeycomb lattice: (a) $w_{50,3}(t)$-$w_{50,3}^{ref}(t)$;  (b) $w_{28,3}(t)$-$w_{28,3}^{ref}(t)$.}\label{f:deviation_FPU}
\end{figure}

\begin{figure}
  \centering
  \includegraphics[width=10cm]{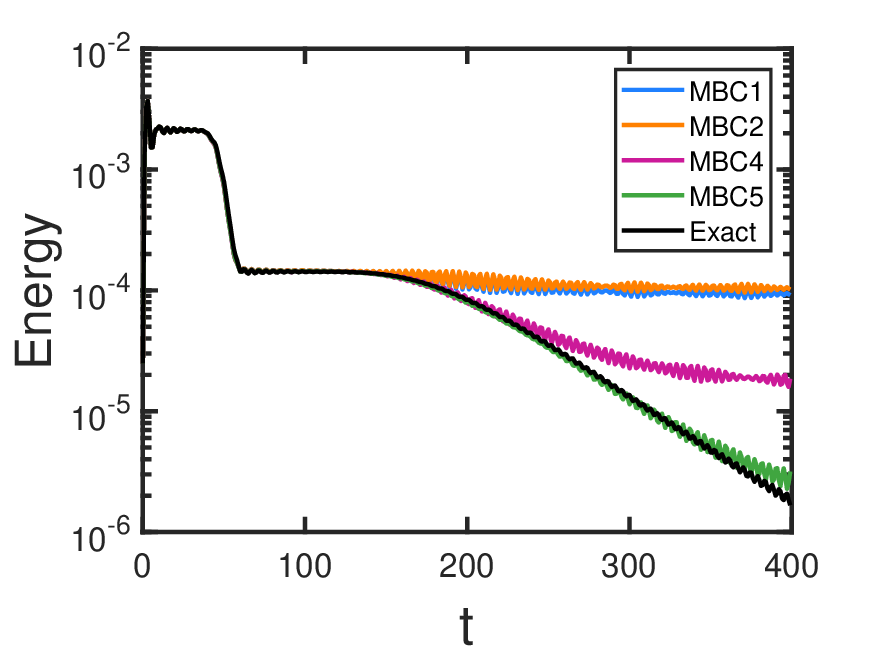}%
  \caption{The kinetic energy of the nonlinear honeycomb lattice during the propagation of the Gaussian hump.}\label{f:energy_FPU}
\end{figure}

We also plot deviations of four matching boundary conditions at two representative atoms. In Fig. \ref{f:deviation_FPU}(a), during long waves passing outward within time $t=0\sim 120$, the curves of four matching boundary conditions almost coincide. That means the low-order MBC and the high-order MBC both well treat long waves at the vertical incident angle. At the same time, for the corner atoms, the amplitude of the pink line (MBC1) and blue line (MBC2) are bigger than the yellow line (MBC5) which is caused by reflections of long waves with large incident angles. After time $t=120$,  the corner reflection of short waves looks smaller than the normal direction. This is because short waves mainly propagate along the normal direction which can be observed in Fig. \ref{f:nonlinear_t200}, so that the reflections at the center of the bottom boundary or top boundary are stronger than other regions. Among the colored lines, the amplitude of the yellow line is the smallest, indicating that the MBC5 provides the best performance for suppressing reflections.

Finally, we plot the kinetic energy curves of the nonlinear honeycomb lattice as shown in Fig. \ref{f:energy_FPU}. The trend  of curves for the nonlinear honeycomb lattice are similar with the harmonic case. In the first platform and the second platform, all colored curves coincide well with the black  one. However, the blue line and the orange line decline very slowly, illustrating that the MBC1 and MBC2 can only eliminate few energy of short waves. The energy for the MBC4 first decays rapidly and then slowly, for the MBC4 is able to deal with part of short waves.  %The purple line first decreases quickly and then the rate of decline become slowly, which means the MBC4 can treat the most of short waves.
The green line is close to the reference solution at all stages which indicates that the MBC5 can well suppress reflections of long waves and short waves simultaneously. All energy curves keeping a downward trend illustrates the stability of the matching boundary conditions for the nonlinear honeycomb lattice.

In addition, remark that the effectiveness of the matching boundary condition is affected by nonlinear strength, for the proposed boundary conditions are designed based on harmonic cases. After many trials, we find that if enhancing the nonlinear effect by increasing the coefficient $\beta$ or the initial amplitude, the boundary reflections will increase correspondingly. However, in practical simulations, the harmonic part plays a leading role in the potential generally, the matching boundary condition can work well to suppress boundary reflections for the nonlinear honeycomb lattice in most applications.

\section{4. Summary}

In this paper, we propose a series of matching boundary conditions for a two-dimensional honeycomb lattice to realize efficient and accurate simulations in a finite domain. The matching boundary condition is written in a linear form composed of the velocities and displacements of several layers of atoms near boundaries with coefficients solved by matching the dispersion relation. An explicit expression makes it convenient for implementations. Only local information involved in the formulas in both space and time guarantees a small cost for the proposed boundary conditions. Furthermore, linear and nonlinear examples verify its effectiveness and stability numerically.

Four representative matching boundary conditions are given in this paper for different cases: the low-order MBC1 and MBC2, and the high-order MBC4 and MBC5. The difference comes from the number of boundary atoms involved and coefficients to match the dispersion relation. We perform a reflection analysis and numerical tests including linear (harmonic potential) and nonlinear (FPU-$\beta$ potential) examples to evaluate these four boundary conditions. These theoretical and numerical analyses show that all the matching boundary conditions can well treat reflections of long waves. However, the low-order MBC1 and MBC2 can not efficiently deal with short waves well especially for the optical branch. In contrast, the high-order MBC4 and MBC5 can efficiently suppress both long waves and short waves simultaneously, because there are more atoms involved and additional coefficients to match the dispersion relation at the optical branch. Among these boundaries, the MBC5 is the most complicated but has the best performance. The trend of the energy curves indicates the stability of the proposed matching boundary conditions. As a result, the low-order MBC1 and MBC2 are simple, and applicable for treating long wave problems with few short waves, such as moving dislocation, or other large deformation in the lattice. The high-order MBC4 and MBC5 have complex forms, but work well in most cases including the optical branch. Thus they can be applied to study the effect of incident waves with high frequency on defects in lattice, or heat transformation problems.

Furthermore, the matching boundary condition is the key component for the heat jet approach to control the system temperature \cite{B. Liu2017, B. Liu2015}. An efficient matching boundary condition can help to design an accurate heat jet approach. Thus the high-order MBC4 and MBC5 proposed in this work may provide a new reliable numerical platform to design heat jet approach for the honeycomb lattice, or honeycomb-like lattice such as graphene, carbon nanotube or singlecrystal silicon in (111) layer. Moreover, in multiscale computations, accurate molecular dynamic simulations rely on non-reflecting artificial boundary conditions. In future work, we desire to propose a multiscale method for the two dimensional honeycomb lattice based on the matching boundary conditions.

\section{Acknowledgments}
B. Liu is supported by the National Natural Science Foundation of China under  Grant No. 12102282. G. Pang is supported by the National Natural Science Foundation of China under Grant No. 11502028. L. Zhang and S. Tang are supported by the National Natural Science Foundation of China Basic Science Center Program for ``Multiscale Problems in Nonlinear Mechanics'' under Grant No. 11988102. L. Zhang is also supported by the National Natural Science Foundation of China under Grant No. 12202451.

\newpage

\section*{Appendix A}
\label{appen:DetailMBC}

As show in Fig. \ref{f:boundary-4}(b), the MBC2 involves two layers of atoms. The MBC2 takes form as
\begin{eqnarray}\nonumber
&&~~~\dot{v}_{n,1}+b_{1,1}(\dot{w}_{n-1,2}+\dot{w}_{n+1,2})
+b_{1,2}(\dot{v}_{n-1,2}+\dot{v}_{n+1,2}) \\
&&=c_{0,0}v_{n,1}+c_{1,1}(w_{n-1,2}+w_{n-1,2})+c_{1,2}(v_{n-1,2}+v_{n-1,2}).
\end{eqnarray}
The three linear equations at long wave limit for the MBC2 are
\begin{numcases}{}\label{eq:mbc2-1}
D_0=c_{0,0}+2c_{1,1}+2c_{1,2}=0,\\\label{eq:mbc2-2}
D_1=3 + 6(b_{1,1}+b_{1,2})-2\sqrt{3}(c_{0,0}+3c_{1,1}+5c_{1,2})= 0, \\\label{eq:mbc2-3}
D_2=2\sqrt{3}+ 2\sqrt{3}(3b_{1,1}+5b_{1,2})-(3c_{0,0}+11c_{1,1}+27c_{1,2})=0.
\end{numcases}
The other two linear equations at wave vector $(0,0.5)$ are
\begin{eqnarray}\label{eq:mbc2-4}
\rm{Re}\{\Delta(0,0.5)\}= 0.7537+ 2.0767b_{1,1} + 2.7069b_{1,2} + 2.2957c_{0,0} + 3.5700c_{1,1} + 0.6535c_{1,2}=0,\\\label{eq:mbc2-5}
\rm{Im}\{\Delta(0,0.5)\}= 1.1357 +1.7662b_{1,1} + 0.3233b_{1,2} - 1.5235c_{0,0} - 4.1977c_{1,1} - 5.4716c_{1,2}=0.
\end{eqnarray}
Solving the five linear equations Eq.(\ref{eq:mbc2-1})-(\ref{eq:mbc2-5}), we obtain the boundary coefficients of the MBC2
\begin{eqnarray}\nonumber
b_{1,1}=0.1671,~ b_{1,2}=0.0139, \\
c_{0,0}=-2.6630, ~c_{1,1}=1.4075, ~c_{1,2}=-0.0760.
\end{eqnarray}

As show in Fig. \ref{f:boundary-4}(c), the MBC4 takes form as
\begin{eqnarray}\nonumber &&~~~\dot{v}_{n,1}+b_{1,1}(\dot{w}_{n-1,2}+\dot{w}_{n+1,2})+b_{1,2}(\dot{v}_{n-1,2}+\dot{v}_{n+1,2}) \\ \nonumber
&&~~~+b_{2,3}(\dot{w}_{n-2,3}+2\dot{w}_{n,3}+\dot{w}_{n+2,3})+b_{2,4}(\dot{v}_{n-2,3}+2\dot{v}_{n,3}+\dot{v}_{n+2,3}) \\ \nonumber
&&=c_{0,0}v_{n,1}+c_{1,1}(w_{n-1,2}+w_{n+1,2})+c_{1,2}(v_{n-1,2}+v_{n+1,2}) \\
&&~~~+c_{2,3}(w_{n-2,3}+2w_{n,3}+w_{n+2,3})+c_{2,4}(v_{n-2,3}+2v_{n,3}+v_{n+2,3}). \label{eq:mbc4-form}
\end{eqnarray}
The three linear equations at long wave limit for the MBC4  are
\begin{numcases}{}\label{eq:mbc4-1}
D_0=c_{0,0}+2c_{1,1}+2c_{1,2}+4c_{2,3}+4c_{2,4}=0,\\\label{eq:mbc4-2}
\nonumber D_1= 3+ 6(b_{1,1} + b_{1,2} + 2b_{2,3} + 2b_{2,4}) \\
~~~~~~~- 2\sqrt{3}(c_{0,0}+3c_{1,1}+5c_{1,2}+12c_{2,3}+16c_{2,4})=0 ,\\\label{eq:mbc4-3}
\nonumber D_2=2\sqrt{3}+2\sqrt{3}(3b_{1,1} + 5b_{1,2}+ 12b_{2,3}+ 16b_{2,4})\\
~~~~~~~-(3c_{0,0}+ 11c_{1,1}+ 27c_{1,2}+ 76c_{2,3}+ 132c_{2,4})=0.
\end{numcases}
The other six linear equations at wave vectors $(0,0.5), (0,1)$ and $(0,2.4)$ are as follows
\begin{eqnarray}\nonumber
\rm{Re}\{\Delta(0,0.5)\}&=& 0.7537+2.0767b_{1,1} + 2.7069b_{1,2} + 5.3815b_{2,3} + 3.9999b_{2,4} \\\label{eq:mbc4-4}
&&+ 2.2957c_{0,0} + 3.5700c_{1,1} + 0.6535c_{1,2} - 1.7695c_{2,3} - 7.4894c_{2,4}= 0,  \\\nonumber
\rm{Im}\{\Delta(0,0.5)\}&=&1.1357+1.7662b_{1,1} + 0.3233b_{1,2} - 0.8754b_{1,2} - 3.7051b_{2,4} \\\label{eq:mbc4-5}
&&- 1.5235c_{0,0} - 4.1977c_{1,1} - 5.4716c_{1,2} - 10.8780c_{2,3} - 8.0853c_{2,4} =0,
\end{eqnarray}
\begin{eqnarray}\nonumber
\rm{Re}\{\Delta(0,1)\}&=& 1.8857+3.9364b_{1,1} + 0.6746b_{1,2} - 2.5281b_{2,3} - 7.9760b_{2,4} \\\label{eq:mbc4-6}
&&+ 0.6789c_{0,0} - 0.6703c_{1,1} - 4.1149c_{1,2} - 7.9196c_{2,3} - 0.0728c_{2,4} = 0, \\\nonumber
\rm{Im}\{\Delta(0,1)\}&=&0.6485- 0.6403b_{1,1} - 3.9307b_{1,2} - 7.5651b_{2,3} - 0.0696b_{2,4}\\\label{eq:mbc4-7}
&&- 1.9741c_{0,0} - 4.1209c_{1,1} - 0.7063c_{1,2} + 2.6465c_{2,3} + 8.3498c_{2,4} = 0,
\end{eqnarray}
\begin{eqnarray}\nonumber
\rm{Re}\{\Delta(0,2.4)\}&=& -3.6841+6.2636b_{1,1} + 4.087b_{1,2} - 13.2119b_{2,3} + 6.1154b_{2,4} \\\label{eq:mbc4-8}
&&- 0.0547c_{0,0} + 1.7932c_{1,1} - 2.8301c_{1,2} + 3.0183c_{2,3} + 6.1882c_{2,4} = 0, \\\nonumber
\rm{Im}\{\Delta(0,2.4)\}&=& -0.1186+3.8875b_{1,1} - 6.1354b_{1,2} + 6.5437b_{2,3} + 13.4161b_{2,4} \\\label{eq:mbc4-9}
&&+ 1.6993c_{0,0} - 2.8891c_{1,1} - 1.8851c_{1,2} + 6.0942c_{2,3} - 2.8208c_{2,4} = 0.
\end{eqnarray}
Solving the nine linear equations Eq.(\ref{eq:mbc4-1})-(\ref{eq:mbc4-9}), we obtain the boundary coefficients of the MBC4
\begin{eqnarray}\nonumber
b_{1,1}=1.6261, ~b_{1,2}=0.5183, ~b_{2,3}=1.0146, ~b_{2,4}=0.1034,\\
c_{0,0}=-4.3763, ~c_{1,1}=0.7035, ~c_{1,2}=0.2887, ~c_{2,3}=0.0729, ~c_{2,4}=0.5250.
\end{eqnarray}

As show in Fig. \ref{f:boundary-4}(d), the MBC5 takes form as
\begin{eqnarray}\nonumber &&~~~\dot{v}_{n,1}+b_{1,1}(\dot{w}_{n-1,2}+\dot{w}_{n+1,2})+b_{1,2}(\dot{v}_{n-1,2}+\dot{v}_{n+1,2}) \\ \nonumber
&&~~~+b_{2,3}(\dot{w}_{n-2,3}+2\dot{w}_{n,3}+\dot{w}_{n+2,3})+b_{2,4}(\dot{v}_{n-2,3}+2\dot{v}_{n,3}+\dot{v}_{n+2,3}) \\ \nonumber
&&~~~+b_{1,5}(\dot{w}_{n-1,4}+\dot{w}_{n+1,4})-b_{3,5}(\dot{w}_{n-3,4}+\dot{w}_{n+3,4})
 \\ \nonumber
&&=c_{0,0}v_{n,1}+c_{1,1}(w_{n-1,2}+w_{n+1,2})+c_{1,2}(v_{n-1,2}+c_{n+1,2}) \\ \nonumber
&&~~~+c_{2,3}(w_{n-2,3}+2w_{n,3}+w_{n+2,3})+c_{2,4}(v_{n-2,3}+2v_{n,3}+v_{n+2,3}) \\
&&~~~+c_{1,5}(w_{n-1,4}+w_{n+1,4})+c_{3,5}(w_{n-3,4}+w_{n+3,4}).
\end{eqnarray}
The three linear equations at long wave limit for the MBC5 are
\begin{numcases}{}\label{eq:mbc5-1}
D_0=c_{0,0}+2c_{1,1}+2c_{1,2}+4c_{2,3}+4c_{2,4}+2c_{1,5}+2c_{3,5}=0,\\\label{eq:mbc5-2}
 \nonumber D_1=3+6(b_{1,1} + b_{1,2} + 2b_{2,3} + 2b_{2,4} + b_{1,5} + b_{3,5} )\\
~~~~~~~ - 2\sqrt{3} (c_{0,0} +3c_{1,1}+5c_{1,2} +12c_{2,3}+16c_{2,4}+9c_{1,5}+9c_{3,5})=0,\\\label{eq:mbc5-3}
 \nonumber D_2=2\sqrt{3}+2\sqrt{3}(3b_{1,1} + 5b_{1,2} + 12b_{2,3} + 16b_{2,4} +9b_{1,5} +9b_{3,5})\\
~~~~~~-(3c_{0,0} + 11c_{1,1} + 27c_{1,2} + 76c_{2,3} + 132c_{2,4} + 83c_{1,5} + 83c_{3,5}) =0.
\end{numcases}
The other ten linear equations at wave vectors $(0,0.5), (0,1), (\sqrt{3}/3,1), (0,2.4)$ and  $(0,3)$ are as follows
\begin{eqnarray} \nonumber\label{eq:mbc5-4}
\rm{Re}\{\Delta(0,0.5)\}&=&0.7537+2.0767b_{1,1} + 2.7069b_{1,2} + 5.3815b_{2,3} + 3.9999b_{2,4} \\ \nonumber
&&+ 1.4098b_{1,5} + 1.4098b_{3,5} + 2.2957c_{0,0} + 3.5700c_{1,1} + 0.6535c_{1,2} \\
&&- 1.7695c_{2,3} - 7.4894c_{2,4} - 4.7164c_{1,5} - 4.7164c_{3,5} = 0, \\ \nonumber\label{eq:mbc5-5}
\rm{Im}\{\Delta(0,0.5)\}&=&1.1357+1.766b_{1,1} + 0.3233b_{1,2} - 0.8754b_{2,3} - 3.7051b_{2,4} \\ \nonumber
&&- 2.3333b_{1,5} - 2.3333b_{3,5} - 1.5235c_{0,0} - 4.1977c_{1,1} - 5.4716c_{1,2} \\
&&- 10.8780c_{2,3} - 8.0853c_{2,4} - 2.8498c_{1,5} - 2.8498c_{3,5} = 0,
\end{eqnarray}
\begin{eqnarray} \nonumber\label{eq:mbc5-6}
\rm{Re}\{\Delta(0,1)\}&=&1.8857+3.9364b_{1,1} + 0.6747b_{1,2} - 2.5281b_{2,3} - 7.9760b_{2,4} \\\nonumber
&&- 3.5305b_{1,5} - 3.5305b_{3,5} + 0.6789c_{0,0} - 0.6703c_{1,1} - 4.1149c_{1,2} \\
&&- 7.9196c_{2,3} - 0.07286c_{2,4} + 1.9419c_{1,5} + 1.9419c_{3,5} = 0, \\\nonumber\label{eq:mbc5-7}
\rm{Im}\{\Delta(0,1)\}&=&0.6485- 0.6403b_{1,1} - 3.9307b_{1,2} - 7.5651b_{2,3} - 0.0696b_{2,4} \\\nonumber
&&+ 1.8549b_{1,5} + 1.8549b_{3,5} - 1.9741c_{0,0} - 4.1209c_{1,1} - 0.7063c_{1,2} \\
&&+ 2.6465c_{2,3} + 8.3498c_{2,4} + 3.6960c_{1,5} + 3.6960c_{3,5} = 0,
\end{eqnarray}
\begin{eqnarray} \nonumber\label{eq:mbc5-8}
\rm{Re}\{\Delta(\tan\frac{\pi}{6},1)\}&=&1.8372+3.4454b_{1,1} + 0.8017b_{1,2} - 1.9418b_{2,3} - 6.1115b_{2,4} \\ \nonumber
&&- 3.0901b_{1,5} - 0.2491b_{3,5} + 0.7182c_{0,0} - 0.5285c_{1,1} - 3.2036c_{1,2} \\
&&- 5.4794c_{2,3} - 0.4069c_{2,4} + 1.5310c_{1,5} + 0.1234c_{3,5} = 0, \\ \nonumber\label{eq:mbc5-9}
\rm{Im}\{\Delta(\tan\frac{\pi}{6},1)\}&=&0.7616- 0.5605b_{1,1} - 3.3974b_{1,2} - 5.8108b_{2,3} - 0.4315b_{2,4} \\ \nonumber
&&+ 1.6236b_{1,5} + 0.1309b_{3,5} - 1.7324c_{0,0} - 3.2489c_{1,1} - 0.7560c_{1,2} \\
&&+ 1.8311c_{2,3} + 5.7629c_{2,4} + 2.9139c_{1,5} + 0.2349c_{3,5} = 0,
\end{eqnarray}
\begin{eqnarray} \nonumber\label{eq:mbc5-10}
\rm{Re}\{\Delta(0,2.4)\}&=&-3.6841+6.2636b_{1,1} + 4.0870b_{1,2} - 13.2123b_{2,3} + 6.1154b_{2,4} \\ \nonumber
&&+ 0.7037b_{1,5} + 0.7037b_{3,5} - 0.0547c_{0,0} + 1.7932c_{1,1} - 2.8299c_{1,2} \\
&&+ 3.0183c_{2,3} + 6.1882c_{2,4} - 3.3848c_{1,5} - 3.3848c_{3,5} = 0, \\ \nonumber\label{eq:mbc5-11}
\rm{Im}\{\Delta(0,2.4)\}&=&-0.1185+3.8875b_{1,1} - 6.1354b_{1,2} + 6.5437b_{2,3} + 13.4159b_{2,4} \\ \nonumber
&&- 7.3383b_{1,5} - 7.3383b_{3,5} + 1.6993c_{0,0} - 2.8891c_{1,1} - 1.8851c_{1,2} \\
&&+ 6.0942c_{2,3} - 2.8208c_{2,4} - 0.3246c_{1,5} - 0.32460c_{3,5} = 0,
\end{eqnarray}
\begin{eqnarray} \nonumber\label{eq:mbc5-12}
\rm{Re}\{\Delta(0,3)\}&=&-4.1970+10.9929b_{1,1} - 12.0052b_{1,2} + 20.4518b_{2,3} - 5.5474b_{2,4} \\ \nonumber
&&- 1.4805b_{1,5} - 1.4805b_{3,5} + 1.9302c_{0,0} - 2.4365c_{1,1} - 1.3391c_{1,2} \\
&&+ 5.9440c_{2,3} - 10.2122c_{2,4} + 5.2011c_{1,5} + 5.2011c_{3,5} = 0,\\ \nonumber\label{eq:mbc5-13}
\rm{Im}\{\Delta(0,3)\}&=&4.5756- 5.7757b_{1,1} - 3.1744b_{1,2} + 14.0902b_{2,3} - 24.2082b_{2,4} \\ \nonumber
&&+ 12.3293b_{1,5} + 12.3293b_{3,5} + 1.7705c_{0,0} - 4.6374c_{1,1} + 5.0644c_{1,2} \\
&&- 8.6276c_{2,3} + 2.3402c_{2,4} + 0.6245c_{1,5} + 0.6245c_{3,5} = 0.
\end{eqnarray}
Solving the thirteen linear equations Eq.(\ref{eq:mbc5-1})-(\ref{eq:mbc5-13}), we obtain the boundary coefficients of the MBC5
\begin{eqnarray}\nonumber
b_{1,1}=4.1088, ~b_{1,2}=5.7272, ~b_{2,3}=2.8636, \\\nonumber
~b_{2,4}=2.0544, ~b_{1,5}=1.0611, ~b_{3,5}=-0.5611,\\\nonumber
c_{0,0}=-6.4564, ~c_{1,1}=-3.2987, ~c_{1,2}=3.2857, \\
~c_{2,3}=-1.6428, ~c_{2,4}=1.6493, ~c_{1,5}=3.8084, ~c_{3,5}=-0.5802.
\end{eqnarray}

\section*{Appendix B}
\label{appen:MBC3}

\begin{figure}
  \centering
  ~~~~~~~~\includegraphics[width=6cm]{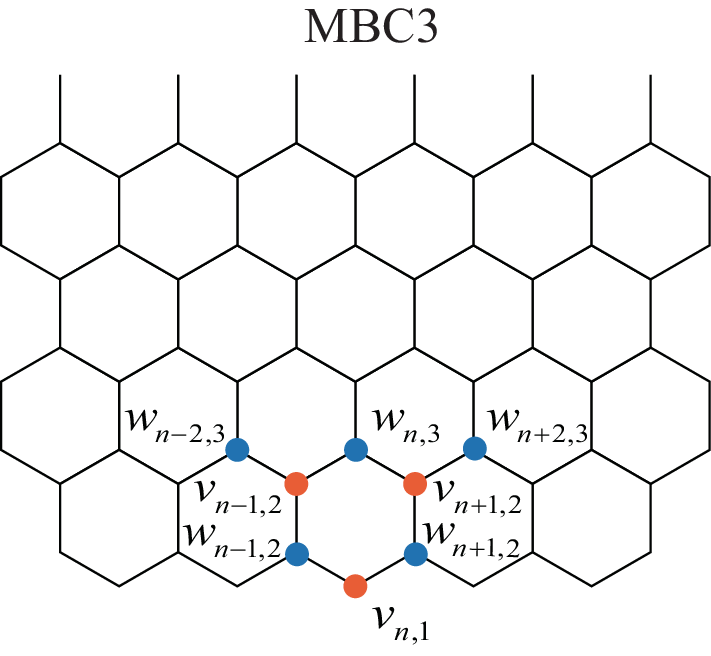} ~~~~~~~~~~~~~~
  \includegraphics[width=8cm]{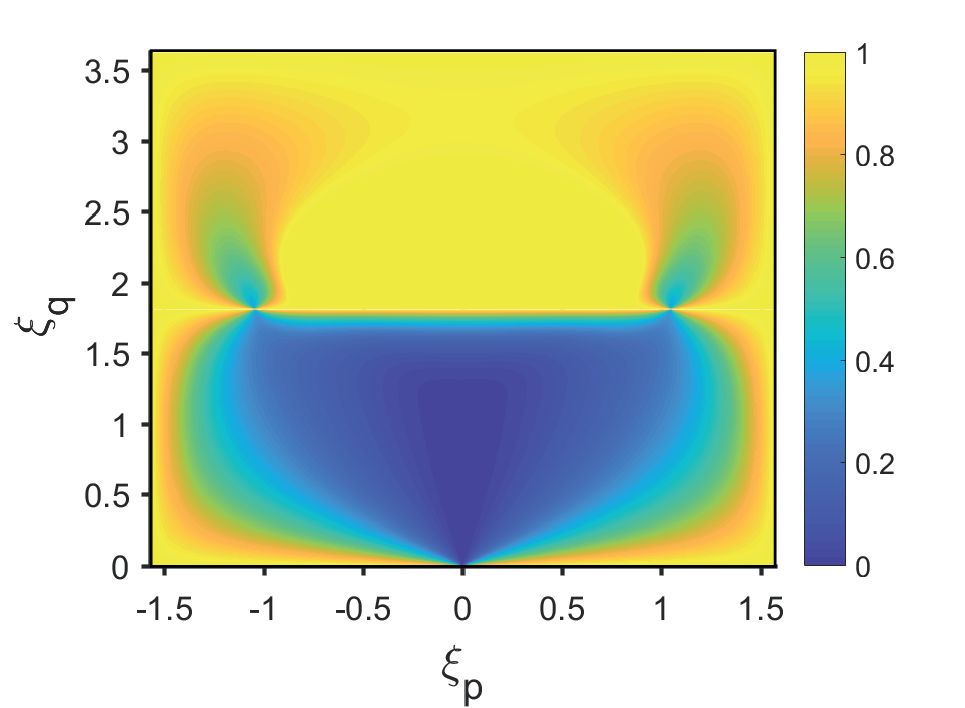}
  \caption{Left: schematic plot of the MBC3 for the bottom boundary; Right: the modulus of reflection coefficients for the MBC3.}\label{f:boundary-mbc3}
\end{figure}
\begin{figure}
  \centering
  \includegraphics[width=8cm]{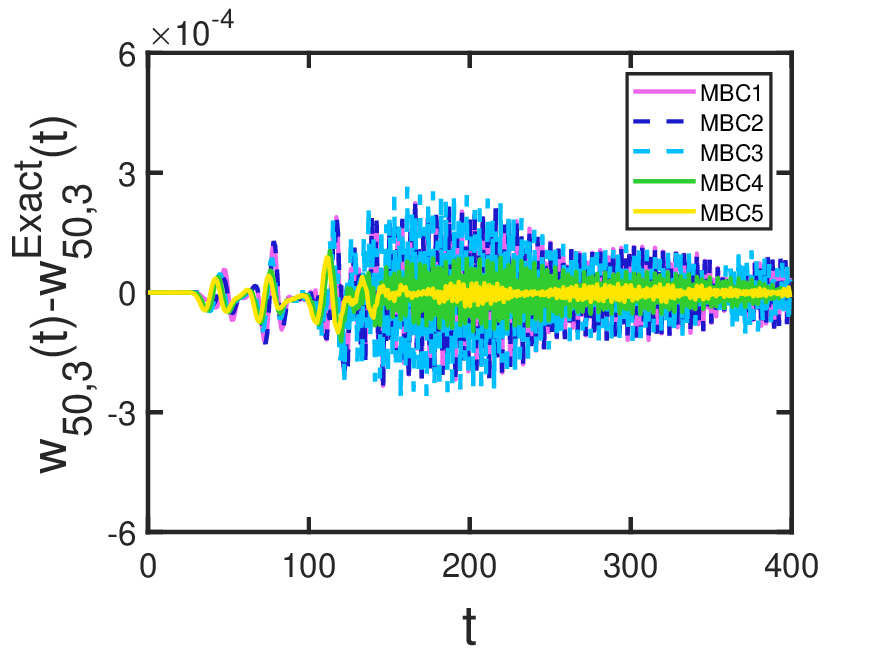} ~~~~~
  \includegraphics[width=8cm]{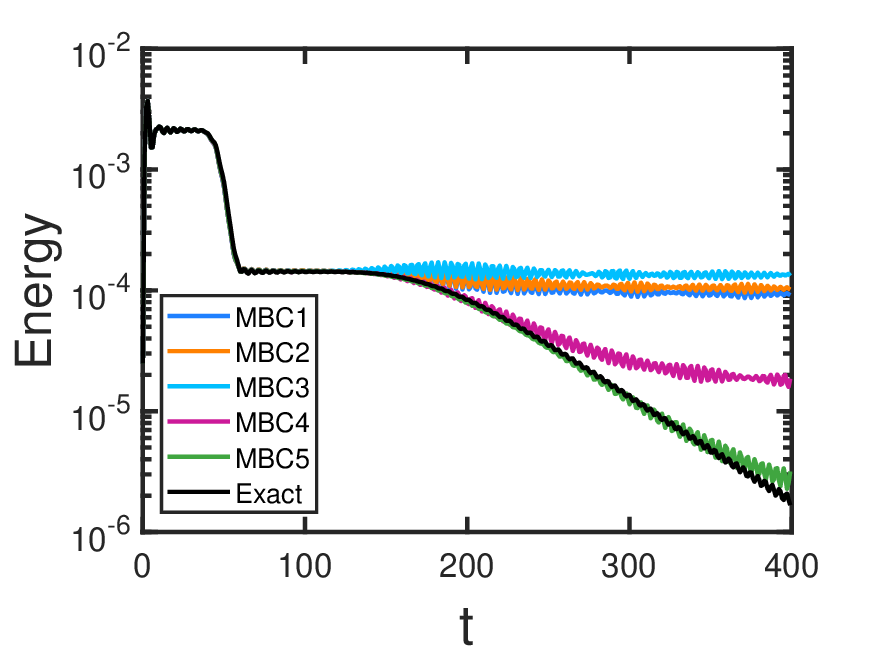}
  \caption{Left: the deviation of displacement $w_{50,3}(t)$-$w_{50,3}^{ref}(t)$ for the harmonic honeycomb lattice with the MBC3; Right: the kinetic energy of the harmonic honeycomb lattice with the MBC3.}\label{f:energy-mbc3}
\end{figure}

The schematic plot of the MBC3 is shown in the left subplot of Fig. \ref{f:boundary-mbc3}. It takes form as
\begin{eqnarray}\nonumber &&~~~\dot{v}_{n,1}+b_{1,1}(\dot{w}_{n-1,2}+\dot{w}_{n+1,2})+b_{1,2}(\dot{v}_{n-1,2}+\dot{v}_{n+1,2})
+b_{2,3}(\dot{w}_{n-2,3}+2\dot{w}_{n,3}+\dot{w}_{n+2,3}) \\ \nonumber
&&=c_{0,0}v_{n,1}+c_{1,1}(w_{n-1,2}+w_{n+1,2})+c_{1,2}(v_{n-1,2}+v_{n+1,2})
+c_{2,3}(w_{n-2,3}+2w_{n,3}+w_{n+2,3}). \label{eq:mbc3}
\end{eqnarray}

The three linear equations at long wave limit for the MBC3 are
\begin{numcases}{}\label{eq:mbc3-1}
D_0=c_{0,0}+2c_{1,1}+2c_{1,2}+4c_{2,3} = 0,\\\label{eq:mbc3-2}
D_1= 3+6(b_{1,1} + b_{1,2} + 2b_{2,3}) - 2\sqrt{3}(c_{0,0}+3c_{1,1}+5c_{1,2}+12c_{2,3}) = 0,\\\label{eq:mbc3-3}
D_2=2\sqrt{3}+2\sqrt{3}(3b_{1,1} +5b_{1,2} +12b_{2,3}) -(3c_{0,0} + 11c_{1,1} + 27c_{1,2} + 76c_{2,3}) = 0.
\end{numcases}
The other four linear equations at wave vectors $(0,0.5), (0,1)$ are as follows
\begin{eqnarray} \nonumber\label{eq:mbc3-4}
\rm{Re}\{\Delta(0,0.5)\}&=& 0.7537+ 2.0767b_{1,1} + 2.7069b_{1,2} + 5.3815b_{2,3} \\
&& + 2.2957c_{0,0}+ 3.5700c_{1,1} + 0.6535c_{1,2} - 1.7695c_{2,3} = 0, \\\nonumber\label{eq:mbc3-5}
\rm{Im}\{\Delta(0,0.5)\}&=&  1.1357 + 1.7661b_{1,1} + 0.3233b_{1,2} - 0.8754b_{2,3} \\
&& - 1.5235c_{0,0}- 4.1977c_{1,1} - 5.4716c_{1,2} - 10.8780c_{2,3} = 0,
\end{eqnarray}
\begin{eqnarray}\nonumber\label{eq:mbc3-6}
\rm{Re}\{\Delta(0,1)\}&=& 1.8857+3.9364b_{1,1} + 0.6746b_{1,2} - 2.5281b_{2,3} \\
&& + 0.6789c_{0,0} - 0.6703c_{1,1} - 4.1149c_{1,2} - 7.9196c_{2,3} = 0, \\\nonumber\label{eq:mbc3-7}
\rm{Im}\{\Delta(0,1)\}&=& 0.6485- 0.6403b_{1,1} - 3.9307b_{1,2} - 7.5651b_{2,3} \\
&& - 1.9741c_{0,0}- 4.1209c_{1,1} - 0.7063c_{1,2} + 2.6465c_{2,3} = 0.
\end{eqnarray}
Solving the seven linear equations Eq.(\ref{eq:mbc3-1})-(\ref{eq:mbc3-7}), we obtain the boundary coefficients of the MBC3
\begin{eqnarray}\nonumber
b_{1,1}=0.8551, ~b_{1,2}=0.8551, ~b_{2,3}=0.2500, \\
c_{0,0}=-3.6755, ~c_{1,1}=1.3285, ~c_{1,2}=-1.3285, ~c_{2,3}=0.9189.
\end{eqnarray}

The modulus of the reflection coefficient is shown in the right subplot of Fig. \ref{f:boundary-mbc3}. Comparing with Fig. \ref{f:reflection}, the blue domain of the MBC3 is larger and darker than the MBC1 and MBC2 for the acoustic branch. However, the color of the optical branch domain is lighter.

Further, see the deviation of displacement in the left subplot of Fig. \ref{f:energy-mbc3} within time $t=0\sim 120$. The deviation of the MBC3 is smaller than the MBC1 and MBC2, while larger than the MBC4 and MBC5. That means the effectiveness of the MBC3 for suppressing long wave reflection is good. However, the MBC3 can not treat short wave reflections. The deviations of the MBC3 become larger than the MBC1 and the MBC2 after $t=120$.

At last, as shown in the right subplot of Fig. \ref{f:energy-mbc3}, the energy of the MBC3 is larger than the MBC1 and MBC2 due to short wave reflections, which coincides with the analysis of the reflection coefficients and the deviation of displacements. We can observe that the MBC3 has no significant advantages compared with the MBC1 and MBC2. Thus, in atomic simulations, it is preferable to use the low-order MBC1 to treat long wave reflections.

\section*{Appendix C}
\label{appen:OtherBoundary}

The MBC4 for other five zigzag boundaries as shown in Fig. \ref{f:computing_domain} are listed as below.

The top boundary is
\begin{eqnarray}\nonumber &&~~~\dot{w}_{n,M}+1.6261(\dot{v}_{n-1,M-1}+\dot{v}_{n+1,M-1})+0.5183(\dot{w}_{n-1,M-1}+\dot{w}_{n+1,M-1}) \\ \nonumber
&&~~~+1.0146(\dot{v}_{n-2,M-2}+2\dot{v}_{n,M-2}+\dot{v}_{n+2,M-2})+0.1034(\dot{w}_{n-2,M-2}+2\dot{w}_{n,M-2}+\dot{w}_{n+2,M-2}) \\ \nonumber
&&=-4.3763 w_{n,M}+0.7035(v_{n-1,M-1}+v_{n+1,M-1})+0.2887(w_{n-1,M-1}+w_{n+1,M-1}) \\
&&~~~+0.0729(v_{n-2,M-2}+2v_{n,M-2}+v_{n+2,M-2})+0.5250(w_{n-2,M-2}+2w_{n,M-2}+w_{n+2,M-2}),
\end{eqnarray}
with
\begin{eqnarray}\nonumber
n=(M+1)/2,~ (M+1)/2+2,~ \cdots, ~N-(M-1)/2-2,~ N-(M-1)/2.
\end{eqnarray}

The lower left boundary is
\begin{eqnarray}\nonumber &&~~~\dot{w}_{n,m}+1.6261(\dot{v}_{n,m}+\dot{v}_{n+1,m-1})+0.5183(\dot{w}_{n+1,m+1}+\dot{w}_{n+2,m}) \\ \nonumber
&&~~~+1.0146(\dot{v}_{n+1,m+1}+2\dot{v}_{n+2,m}+\dot{v}_{n+3,m-1})+0.1034(\dot{w}_{n+2,m+2}+2\dot{w}_{n+3,m+1}+\dot{w}_{n+4,m}) \\ \nonumber
&&=-4.3763 w_{n,m}+0.7035(v_{n,m}+v_{n+1,m-1})+0.2887(w_{n+1,m+1}+w_{n+2,m}) \\
&&~~~+0.0729(v_{n+1,m+1}+2v_{n+2,m}+v_{n+3,m-1})+0.5250(w_{n+2,m+2}+2w_{n+3,m+1}+w_{n+4,m}),
\end{eqnarray}
with
\begin{eqnarray}\nonumber
&&m=3, ~4, ~\cdots, ~(M+1)/2-2, ~(M+1)/2-1, \\\nonumber
&&n(m)=(M+3)/2-m.
\end{eqnarray}

The lower right boundary is
\begin{eqnarray}\nonumber &&~~~\dot{w}_{n,m}+1.6261(\dot{v}_{n,m}+\dot{v}_{n-1,m-1})+0.5183(\dot{w}_{n-1,m+1}+\dot{w}_{n-2,m}) \\ \nonumber
&&~~~+1.0146(\dot{v}_{n-1,m+1}+2\dot{v}_{n-2,m}+\dot{v}_{n-3,m-1})+0.1034(\dot{w}_{n-2,m+2}+2\dot{w}_{n-3,m+1}+\dot{w}_{n-4,m}) \\ \nonumber
&&=-4.3763 w_{n,m}+0.7035(v_{n,m}+v_{n-1,m-1})+0.2887(w_{n-1,m+1}+w_{n-2,m}) \\
&&~~~+0.0729(v_{n-1,m+1}+2v_{n-2,m}+v_{n-3,m-1})+0.5250(w_{n-2,m+2}+2w_{n-3,m+1}+w_{n-4,m}),
\end{eqnarray}
with
\begin{eqnarray}\nonumber
&&m=3, ~4, ~\cdots, ~(M+1)/2-2, ~(M+1)/2-1, \\\nonumber
&&n(m)=-(M+1)/2+N+m.
\end{eqnarray}

The upper left boundary is
\begin{eqnarray}\nonumber &&~~~\dot{v}_{n,m}+1.6261(\dot{w}_{n+1,m+1}+\dot{w}_{n,m})+0.5183(\dot{v}_{n+2,m}+\dot{v}_{n+1,m-1}) \\ \nonumber
&&~~~+1.0146(\dot{w}_{n+3,m+1}+2\dot{w}_{n+2,m}+\dot{w}_{n+1,m-1})+0.1034(\dot{v}_{n+4,m}+2\dot{v}_{n+3,m-1}+\dot{v}_{n+2,m-2}) \\ \nonumber
&&=-4.3763 v_{n,m}+0.7035(w_{n+1,m+1}+w_{n,m})+0.2887(v_{n+2,m}+v_{n+1,m-1}) \\
&&~~~+0.0729(w_{n+3,m+1}+2w_{n+2,m}+w_{n+1,m-1})+0.5250(v_{n+4,m}+2v_{n+3,m-1}+v_{n+2,m-2}),
\end{eqnarray}
with
\begin{eqnarray}\nonumber
&&m=(M+1)/2+1, ~(M+1)/2+2,~ \cdots, ~M-3, ~M-2, \\\nonumber
&&n(m)=-(M-1)/2+m.
\end{eqnarray}

The upper right boundary is
\begin{eqnarray}\nonumber &&~~~\dot{v}_{n,m}+1.6261(\dot{w}_{n-1,m+1}+\dot{w}_{n,m})+0.5183(\dot{v}_{n-2,m}+\dot{v}_{n-1,m-1}) \\ \nonumber
&&~~~+1.0146(\dot{w}_{n-3,m+1}+2\dot{w}_{n-2,m}+\dot{w}_{n-1,m-1})+0.1034(\dot{v}_{n-4,m}+2\dot{v}_{n-3,m-1}+\dot{v}_{n-2,m-2}) \\ \nonumber
&&=-4.3763 v_{n,m}+0.7035(w_{n-1,m+1}+w_{n,m})+0.2887(v_{n-2,m}+v_{n-1,m-1}) \\
&&~~~+0.0729(w_{n-3,m+1}+2w_{n-2,m}+w_{n-1,m-1})+0.5250(v_{n-4,m}+2v_{n-3,m-1}+v_{n-2,m-2}),
\end{eqnarray}
with
\begin{eqnarray}\nonumber
&&m=(M+1)/2+1, ~(M+1)/2+2,~ \cdots, ~M-3, ~M-2, \\\nonumber
&&n(m)=(M+1)/2+N-m.
\end{eqnarray}
Here, the subscript $n$ of the boundary atoms are the function of $m$ for the oblique boundaries of the hexagon computational domain.  The other matching boundary conditions and the MBC1 for the corner atoms can be written out in the same manner.

\newpage

%
%------------------------------------------------------------------

%\clearpage
%\bibliographystyle{apsrev}
%
%\bibliography{neutronshellgap} % Produces the bibliography via BibTeX.

\begin{thebibliography}{99}


\bibitem{S. Hajilar2015} S. Hajilar and B. Shafei, Nano-scale investigation of elastic properties of hydrated cement paste constituents using molecular dynamics simulations, Comput. Mater. Sci. 101 (2015) 216-226.

\bibitem{A. Kumar2020} A. Kumar, K. Sharma and A. Dixit, A review on the mechanical and thermal properties of graphene and graphene-based polymer nanocomposites: understanding of modelling and MD simulation, Mol. Simulat. 46(2) (2020) 136-154.



\bibitem{S. Groh2009} S. Groh, E. Marin, M. Horstemeyer, et al., Dislocation motion in magnesium: a study by molecular statics and molecular dynamics, Model. Simul. Mater. Sci. Eng. 17(7) (2009) 075009.

\bibitem{E. Oren2016} E. Oren, E. Yahel and G. Makov, Dislocation kinematics: a molecular dynamics study in Cu, Model. Simul. Mater. Sci. Eng. 25(2) (2016) 025002.



\bibitem{Y. Lee2011} Y. Lee, S. Lee and G. Hwang, Effects of vacancy defects on thermal conductivity in crystalline silicon: A nonequilibrium molecular dynamics study, Phys. Rev. B, 83(12) (2011) 125202.

\bibitem{G. Berdiyorov} G. Berdiyorov and  F. Peeters, Influence of vacancy defects on the thermal stability of silicene: a reactive molecular dynamics study, Rsc adv. 4(3) (2014) 1133-1137.

\bibitem{M. Noshin} M. Noshin, A. Khan, I. Navid, et al., Impact of vacancies on the thermal conductivity of graphene nanoribbons: A molecular dynamics simulation study, Aip Adv. 7(1) (2017) 015112.



\bibitem{C. Rountree2002} C. Rountree, R. Kalia, E. Lidorikis, et al., Atomistic aspects of crack propagation in brittle materials: Multimillion atom molecular dynamics simulations, Annual Review of Materials Research, 32(1) (2002) 377-400.

\bibitem{P. Budarapu2015} P. Budarapu, B. Javvaji, V. Sutrakar, et al., Crack propagation in graphene, J Appl. Phys. 118(6) (2015) 064307.

%Dehaghani M Z, Mashhadzadeh A H, Salmankhani A, et al. Fracture toughness and crack propagation behavior of nanoscale beryllium oxide graphene-like structures: a molecular dynamics simulation analysis[J]. Engineering fracture Mechanics, 2020, 235: 107194.

\bibitem{Kadowaki2005} H. Kadowaki and W. Liu, A multiscale approach for the micropolar continuum model, Laser Phys. 15(3) (2005) 269-282.

\bibitem{R. Miller2009} R. Miller and E. Tadmor, A unified framework and performance benchmark of fourteen multiscale atomistic/continuum coupling methods, Model. simul. mater. sci. and eng. 17(5) (2009) 053001.

\bibitem{V. Yamakov2014} V. Yamakov, D. Warner, R. Zamora, et al., Investigation of crack tip dislocation emission in aluminum using multiscale molecular dynamics simulation and continuum modeling, J. Mech. Phys. Solids, 65 (2014)  35-53.

\bibitem{S. Tang2006} S. Tang, T. Hou and W. Liu, A pseudo-spectral multiscale method: interfacial conditions and coarse grid equations, J. Comput. Phys. 213(1) (2006) 57-85.

\bibitem{Adelman&Doll 1974} S. Adelman and J. Doll, Generalized Langevin equation approach for atom/solid-surface scattering: collinear atom/harmonic chain model, J. Chem. Phys. 61(10) (1974) 4242-4245.

\bibitem{Adelman&Doll 1976} S. Adelman and J. Doll, Generalized Langevin equation approach for atom/solid-surface scattering: general formulation for classical scattering off harmonic solids, J. Chem. Phys. 64(6) (1976) 2375-2388.

\bibitem{Berenger1994} J. Berenger, A perfectly matched layer for the absorption of electromagnetic waves, J. comput. phys. 114(2) (1994) 185-200.

\bibitem{W. E2001} W. E and Z. Huang, Matching conditions in atomistic-continuum modeling of materials, Phys. Rev. Lett. 87(13) (2001) 135501.

\bibitem{W. E2002} W. E and Z. Huang, A dynamic atomistic–continuum method for the simulation of crystalline materials, J. Comput. Phys. 182(1) (2002) 234-261.

\bibitem{S. Tang2008} S. Tang, A finite difference approach with velocity interfacial conditions for multiscale computations of crystalline solids, J. Comput. Phys. 227(8) (2008) 4038-4062.

\bibitem{G. Pang ALEX 2012} G. Pang, L. Bian and S. Tang, Almost exact boundary condition for one-dimensional Schr\"{o}dinger equations, Phys. Rev. E, 86(6) (2012) 066709.

\bibitem{X. Wang2010} X. Wang, S. Tang, Matching boundary conditions for diatomic chains, Comput. Mech. 46 (2010) 813-826.

\bibitem{X. Wang2013} X. Wang, S. Tang, Matching boundary conditions for lattice dynamics, Int. J. Numer. Methods Eng. 93(12) (2013) 1255-1285.



\bibitem{W. Liu2000} W. Liu, S. Hao, T. Belytschko, S. Li and C. T. Chang, Multi-scale methods, International Journal for Numerical Methods in Engineering 47(7), 1343-1361 (2000).

\bibitem{W. Liu2003} G. Wagner and W. Liu, Coupling of atomistic and continuum simulations using a bridging scale decomposition, Journal of Computational Physics 190(1) (2003) 249-274.

\bibitem{W. Liu2004} W. Liu, E. Karpov, S. Zhang and H. Park, An introduction to computational nanomechanics and materials, Comput. Method. Appl. Mech. and Eng. 193(17) (2004) 1529-1578.

\bibitem{M. Dreher2008} M. Dreher and S. Tang, Time history interfacial conditions in multiscale computations of lattice oscillations, Comput. Mech. 41 (2008) 683-698.

\bibitem{G. Pang ALEX 2011} G. Pang and S. Tang, Time history kernel functions for square lattice. Comput. Mech. 48 (2011) 699-711.

%\bibitem{Cai2000} W. Cai, M. de Koning, V. V. Bulatov and S. Yip, Minimizing boundary reflections in coupled-domain simulations, Phys. Rev. Lett. 85(15) (2000) 3213.

\bibitem{W. Chew1996} W. Chew and Q. Liu, Perfectly matched layers for elastodynamics: a new absorbing boundary condition, J. Comput. Acoust. 4(04) (1996) 341-359.

\bibitem{Hastings1996} F. D. Hastings, J. B. Schneider and S. L. Broschat, Application of the perfectly matched layer (PML) absorbing boundary condition to elastic wave propagation, J. Acoust. Soc. Am. 100 (5) (1996) 3061-3069.

\bibitem{C. Zheng2007} C. Zheng, A perfectly matched layer approach to the nonlinear Schr\"{o}dinger wave equations. J. Comput. Phys. 227(1) (2007) 537-556.

\bibitem{S. Li2005} A. C. To and S. Li, Perfectly matched multiscale simulations, Phys. Rev. B 72(3) (2005) 035414.

\bibitem{S. Li2006} S. Li, X. Liu, A. Agrawal and A. C. To, Perfectly matched multiscale simulations for discrete lattice systems: extension to multiple dimensions, Phys. Rev. B 74(4) (2006) 045418.

\bibitem{X. Li2006} X. Li and W. E, Variational boundary conditions for molecular dynamics simulations of solids at low temperature, Commun. Comput. Phys. 1(1) (2006) 135-175.

\bibitem{X. Li2007} X. Li and W. E, Variational boundary conditions for molecular dynamics simulations of crystalline solids at finite temperature: treatment of the thermal bath. Phys. Rev. B 76(10) (2007) 104107.

\bibitem{G. Pang ALEX 2017} G. Pang, Y. Yang and S. Tang, Exact boundary condition for semi-discretized Schr\"{o}dinger equation and heat equation in a rectangular domain, J. Sci. Comput. 72 (2017) 1-13.

\bibitem{G. Pang Bessel 2017} G. Pang and S. Tang, Approximate linear relations for Bessel functions. Commun. Math. Sci. 15(7) (2017) 1967-1986.

\bibitem{B. Liu2017} B. Liu, S. Tang and J. Chen, Heat jet approach for finite temperature atomic simulations of triangular lattice, Comput. Mech. 59(5) (2017) 843-859.

\bibitem{M. Fang2012} M. Fang, S. Tang, Z. Li, et al., Artificial boundary conditions for atomic simulations of face-centered-cubic lattice, Comput. Mech. 50 (2012) 645-655.

\bibitem{M. Fang2013} M. Fang, S. Tang, Z. Li, et al. Matching boundary conditions for scalar waves in body-centered-cubic lattices, Adv. Appl. Math. Mech. 5(3) (2013) 337-350.

\bibitem{S. Ji2014} S. Tang and S. Ji, Stability of atomic simulations with matching boundary conditions, Adv. Appl. Math. Mech. 6(5) (2014) 539-551.

\bibitem{B. Liu2015} S. Tang and B. Liu, Heat jet approach for atomic simulations at finite temperature, Commun. Comput. Phys. 18(5) (2015) 1445-1460.

\bibitem{B. Liu2020} B. Liu and S. Tang, Non-equilibrium atomic simulation for Frenkel-Kontorova model with moving dislocation at finite temperature, Chin. Phys. B 29(11) (2020) 110501.

\bibitem{L. Zhang2022} L. Zhang, S. Tang and B. Liu, Finite difference approach for multiscale computations of atomic chain at finite temperature, Comput. Math. Appl. 110 (2022) 77-90.


\bibitem{G. Wagner2004} G. Wagner, E. Karpov and W. Liu, Molecular dynamics boundary conditions for regular crystal lattices, Comput. Meth. Appl. Mech. Eng. 193(17-20) (2004) 1579-1601.

\bibitem{S. Medyanik2006}S. Medyanik, E. Karpov and W. Liu, Domain reduction method for atomistic simulations, J. Comput. Phys. 218(2) (2006) 836-859.

\bibitem{X. Yu2012} X. Yu and C. Sarris, A perfectly matched layer for subcell FDTD and applications to the modeling of graphene structures, IEEE Antenn. Wirel. Propag. Lett. 11 (2012) 1080-1083.

\bibitem{S. Ji2017} S. Ji and S. Tang, Artificial boundary conditions for out-of-plane motion in penta-graphene, Acta Mech. Sin. 33(6) (2017) 992-998.


%\bibitem{G. Pang2021} Pang G, Yang Y, Antoine X, et al. Stability and convergence analysis of artificial boundary conditions for the Schr\"{o}dinger equation on a rectangular domain, Math. Comput. 90(332) (2021) 2731-2756.

\bibitem{G. Pang2018 corner} Pang G, Ji S, Yang Y, et al. Eliminating corner effects in square lattice simulation, Comput. Mech. 62 (2018) 111-122.

\bibitem{S. Lepri2003} S. Lepri, R. Livi and A. Politi, Thermal conduction in classical low-dimensional lattices, Phys. rep. 377(1) (2003) 1-80.

\end{thebibliography}

%\end{CJK*}
\end{document}